\documentclass[a4paper]{article}
    \usepackage[T1]{fontenc}
    \usepackage{mathtools}
    \usepackage{amssymb}
    \usepackage{amsthm}
    \usepackage{authblk}
    \usepackage{algorithm}
    \usepackage{algpseudocodex}
    \usepackage{mathtools}
    \usepackage{bm}
    \usepackage{ulem}
    \usepackage{bbding}
    \usepackage{xfrac}
    \usepackage{subcaption}
    \usepackage{pgfplots}
    \usepackage{titlesec}
    \usepackage{rotating}
    \usepackage{comment}
    \usepackage{cite}
    \usepackage[final]{showlabels}
    \usepackage[textsize=footnotesize,textwidth=2cm]{todonotes}
    \usepackage[hidelinks,hypertexnames=false,bookmarks=false,draft]{hyperref}
    \usepackage[capitalise]{cleveref}
    \usepackage{xcolor}
    \usepackage{xargs}
    \usepackage{colortbl}
    \usepackage{multicol}
    \usepackage{enumerate}
    
    \usepackage{xr}
    \newcommand{\defeq}{\vcentcolon=}
    
    \newcommand*{\tran}{^{\mkern-1.5mu\mathsf{T}}}
    \newcommandx{\axtd}[2][1=]{\todo[inline,linecolor=orange,backgroundcolor=orange!40,bordercolor=orange,#1]{#2}}
    \newcommandx{\dmtd}[2][1=]{\todo[linecolor=magenta,backgroundcolor=magenta!40,bordercolor=magenta,#1]{#2}}
    \newcommand{\cgc}[1]{\cellcolor{green!#1}\Checkmark}
    \newcommand{\cyc}[1]{\cellcolor{yellow!#1}\Checkmark}
    \newcommand{\cyz}[1]{\cellcolor{yellow!#1}$\bm{0}$}
    \newcommand{\crc}[1]{\cellcolor{red!#1}\Checkmark}
    \newcommand{\crz}[1]{\cellcolor{red!#1}$\bm{0}$}
    \newcommand{\cbc}[1]{\cellcolor{blue!#1}\Checkmark}
    \newcommand{\cbz}[1]{\cellcolor{blue!#1}$\bm{0}$}

    \title{Asymptotic analysis of parameterised univariate Gaussian splitting}
    \author{Dmitry Mikhin and Athena Xiourouppa}
    \date{\today}

    \pgfplotsset{compat=1.18}

\makeatletter
\def\user@resume{resume}
\def\user@intermezzo{intermezzo}
\newcounter{previousequation}
\newcounter{lastsubequation}
\newcounter{savedparentequation}
\renewenvironment{subequations}[1][]{%
      \def\user@decides{#1}%
      \setcounter{previousequation}{\value{equation}}%
      \ifx\user@decides\user@resume
           \setcounter{equation}{\value{savedparentequation}}%
      \else
      \ifx\user@decides\user@intermezzo
           \refstepcounter{equation}%
      \else
           \setcounter{lastsubequation}{0}%
           \refstepcounter{equation}%
      \fi\fi
      \protected@edef\theHparentequation{%
          \@ifundefined {theHequation}\theequation \theHequation}%
      \protected@edef\theparentequation{\theequation}%
      \setcounter{parentequation}{\value{equation}}%
      \ifx\user@decides\user@resume
           \setcounter{equation}{\value{lastsubequation}}%
         \else
           \setcounter{equation}{0}%
      \fi
      \def\theequation  {\theparentequation  \alph{equation}}%
      \def\theHequation {\theHparentequation \alph{equation}}%
      \ignorespaces
}{%
  \ifx\user@decides\user@resume
       \setcounter{lastsubequation}{\value{equation}}%
       \setcounter{equation}{\value{previousequation}}%
  \else
  \ifx\user@decides\user@intermezzo
       \setcounter{equation}{\value{parentequation}}%
  \else
       \setcounter{lastsubequation}{\value{equation}}%
       \setcounter{savedparentequation}{\value{parentequation}}%
       \setcounter{equation}{\value{parentequation}}%
  \fi\fi
  \ignorespacesafterend
}
\makeatother

\titleformat{\paragraph}
{\normalfont\normalsize\bfseries}{\theparagraph}{1em}{}
\titlespacing*{\paragraph}
{0pt}{3.25ex plus 1ex minus .2ex}{1.5ex plus .2ex}

\titleformat{\subparagraph}
{\normalfont\normalsize\bfseries}{\thesubparagraph}{1em}{}
\titlespacing*{\subparagraph}
{0pt}{3.25ex plus 1.5ex minus .2ex}{2.5ex plus .2ex}

\makeatletter
\@addtoreset{section}{part}
\makeatother

\titleformat{\subparagraph}
{\normalfont\normalsize\bfseries}{\thesubparagraph}{1em}{}
\titlespacing*{\subparagraph}
{0pt}{3.25ex plus 1ex minus .2ex}{1.5ex plus .2ex}

\makeatletter
\newcommand{\raisemath}[1]{\mathpalette{\raisem@th{#1}}}
\newcommand{\raisem@th}[3]{\raisebox{#1}{$#2#3$}}
\makeatother

\begin{document}

\maketitle

\section{Introduction\label{sec:Intro}}

This document provides in-depth details for the derivation of the univariate splitting
algorithm developed in \cite{Mikhin:2026:Splitting};
the minimal set of definitions is repeated in \cref{tbl:notation}.
The algorithm approximates the standard, $1$-D Gaussian distribution with a mixture of
uniformly spaced homoscedastic Gaussian components. The solution is found by minimising the squared $L^2$ norm
of the mismatch between the approximation and the original Gaussian.

This text presents asymptotic analyses of the proposed
splitting in the limit of small step $h$ between the mixand means and in the limit of large
number of mixands $M$.
These results were originally obtained in the PhD thesis \cite{Xiourouppa:2027:PhD}; this is an
abridged version adapted to work as extended Appendices of \cite{Mikhin:2026:Splitting}. For this
reason, all the section and equation references internal to this document are prefixed with ``A'' or
``B'', as it would be in Appendix A or B,
for example, \cref{eq:dddh:h0:Odd} for the very first equation below. On the contrary, all
references not starting with a letter refer to the main paper \cite{Mikhin:2026:Splitting}, e.g.,
\cref{G1:eq:cd:Odd} refers to the corresponding equation in that text.

\begin{table}[ht]
    \caption{Notation and abbreviations}
    \label{tbl:notation}
    \begin{tabular}{|p{0.18\textwidth}|p{0.75\textwidth}|}
    \hline
        \underline{Abbreviations:} & \\
        lhs      & left-hand side \\
        rhs      & right-hand side \\
        GS       & Gaussian Sum \\
        \underline{Notation:} & \\
        $\sigma$ & the mixand standard deviation; given to the algorithm by the user \\
        $M$      & half-number of mixands; the total number is $2M + 1$ for the odd case and $2M$
                   for the even; given by the user \\
        $L^2$    & the approximation error (Euclidian norm squared) defined as the integral of the
                   squared difference between the standard normal function and its GS approximation,
                   \cref{G1:eq:L2} \\
        $h$      & the step between the mixand means; found by minimisation of $L^2$ \\
        $\bm{w}$ & the length-$M$ vector of mixand weights; found by minimisation of $L^2$; all
                   weights must be positive and their sum equals to one. \\
    \hline
    \end{tabular}
\end{table}

\section{Derivation of \texorpdfstring{$L^2$}{L²} and \texorpdfstring{$\bm{w}$}{w} for small
\texorpdfstring{$h$}{h}}\label{sec:SmallH}

\renewcommand{\theequation}{A.\arabic{equation}}
\setcounter{equation}{0}

This section presents asymptotic analysis of the proposed splitting method in the limit of small
step $h$ between the mixand means. The summary of these results is given in
\cite[\cref{G1:sec:SmallH:Odd,G1:sec:SmallH:Even}]{Mikhin:2026:Splitting}.

\subsection{Odd case}\label{sec:SmallH:Odd}

In our splitting method, the mixand weights $\bm{w}$ are found by solving the linear system of
equations $\bm{A} \bm{w} = \bm{b}$, where the elements of the matrix $\bm{A}$ and the rhs vector
$\bm{b}$ are given by \cref{G1:eq:cd:Odd,G1:eq:Gprods}. At $h = 0$, all functions $d_\alpha$,
$c_{\alpha,\beta}$ in \cref{G1:eq:Gprods} become identical, and therefore, the elements $a_{m,k}$,
$b_m$ in \cref{G1:eq:cd:Odd} turn to zero, rendering the weights undefined.

To obtain the solution for small but finite $h$,
we decompose $\bm{A}$ and $\bm{b}$ into Maclaurin series over $h_* \defeq h / \sigma$.
The derivatives of $d_{\alpha}$ and
$c_{\alpha, \beta}$ from \cref{G1:eq:Gprods} with respect to $h_*$ at $h_* = 0$ are
\begin{align}
    \label{eq:dddh:h0:Odd}
    \frac{\partial^{2n} d_{\alpha}}{\partial h_*^{2n}} \Bigr|_{h_*=0} &= \frac{(-1)^n (2n)! \alpha^{2n} \sigma^{2n}}{n! \, (2(1 + \sigma^2))^{n + 1/2} \sqrt{\pi} }
    \text{,}
    \\
    \label{eq:dcdh:h0:Odd}
    \frac{\partial^{2n} c_{\alpha, \beta}}{\partial h_*^{2n}} \Bigr|_{h_*=0} &= \frac{(-1)^n (2n)! (\alpha - \beta)^{2n}}{n! \, 2^{2n + 1} \sigma \sqrt{\pi}} \text{.}
\end{align}
By symmetry, all odd derivatives at $h_* = 0$ are zero. Substituting \cref{eq:dddh:h0:Odd,eq:dcdh:h0:Odd} into
\cref{G1:eq:cd:Odd}, we obtain for the derivatives of $b_m$ and $a_{m,k}$:
\begingroup
\allowdisplaybreaks
\begin{alignat}{2}
    \label{eq:dbdh:h0:Odd}
        &\frac{1}{(2n)!} \frac{\partial^{2n} b_m}{\partial h_*^{2n}} \Bigr|_{h_*=0} &&=
            \frac{(-1)^n}{n! \, 2^{2n} \sigma \sqrt{\pi}} \left( \frac{2^{n + 1/2} \sigma^{2n + 1}}{(1 + \sigma^2)^{n + 1/2}} - 1 \right) m^{2n} \text{,} \\
    \notag
        &\frac{1}{(2n)!} \frac{\partial^{2n} a_{m,k}}{\partial h_*^{2n}} \Bigr|_{h_*=0} &&=
            \frac{
                (-1)^n \left[ (k - m)^{2n} + (k + m)^{2n} - 2 k^{2n} - 2 m^{2n} \right]
            }{n! \, 2^{2n} \sigma \sqrt{\pi}} \\
    \label{eq:dadh:h0:Odd}
        & &&=
            \frac{(-1)^n}{n! \, 2^{2n - 1} \sigma \sqrt{\pi}} \sum_{j=1}^{n-1} \binom{2n}{2j} k^{2j} m^{2(n - j)} \text{.}
\end{alignat}
\endgroup
The second order derivative of $a_{m,k}$ ($n = 1$ in \cref{eq:dadh:h0:Odd}) is zero. We moved
the factors of $(2n)!$ to the lhs because this is how they would appear in the Maclaurin series.

Using \cref{eq:dbdh:h0:Odd}, we write the decomposition of the rhs vector $\bm{b}$ in compact vector
form as (reproducing \cref{G1:eq:decb:h0:Odd} here)
\begin{equation}
    \label{eq:decb:h0:Odd}
    \bm{b} = \sum_{n=1}^{+\infty} \beta_{2n} \bm{p}_{2n} h_*^{2n} \text{,}
\end{equation}
where
\begin{equation}
    \label{eq:elmb:h0:Odd}
        \bm{b}_{2n} = \beta_{2n} \bm{p}_{2n} \text{, $\forall n \ge 1$,}
\end{equation}
\begin{equation}
    \label{eq:beta:h0:Odd}
         \beta_{2n} = \frac{(-1)^n}{n! \, 2^{2n} \sigma \sqrt{\pi}} \left( \frac{2^{n + 1/2} \sigma^{2n + 1}}{(1 + \sigma^2)^{n + 1/2}} - 1 \right)
         \! \text{,}
\end{equation}
and
\begin{equation}
    \label{eq:p2n:Odd}
    \bm{p}_{2n} = \bigl[ 1, 2^{2n}, \dots, M^{2n} \bigr] \text{.}
\end{equation}
The vectors $\bm{p}_{2n}$ for $n \le M$ are independent, although not mutually orthogonal.

Similarly, using \cref{eq:dadh:h0:Odd} we decompose the matrix $\bm{A}$ as (repeating
\cref{G1:eq:deca:h0:Odd})
\begin{equation}
    \label{eq:deca:h0:Odd}
        \bm{A} = \sum_{n=2}^{+\infty} \bm{A}_{2n} h_*^{2n} \text{,}
\end{equation}
where
\begin{equation}
    \label{eq:A:h0:Odd}
        \bm{A}_{2n} =
            \sum_{j = 1}^{n - 1} \alpha^{(2 n)}_{2j,2(n - j)}
                \left( \bm{p}_{2 j} \otimes \bm{p}_{2(n - j)} \right) \text{, $\forall n \ge 2$,}
\end{equation}
the symbol $\otimes$ denotes the outer product of vectors, and
\begin{equation}
    \label{eq:alpha:h0:Odd}
    \alpha^{(2n)}_{2j, 2(n - j)} = \frac{(-1)^n}{n! \, 2^{2n-1} \sigma \sqrt{\pi}} \binom{2n}{2j} \text{.}
\end{equation}
The second of the lower indices in $\alpha^{(2n)}_{2j,2(n - j)}$ is technically unnecessary, as it
is fully defined by the first lower and the upper indices. However, it makes for convenient notation
relating the indices of $\alpha$ with the indices of $\bm{p}$ in the outer product, for example,
\begin{equation}
    \notag
        \bm{A}_8 =
            \alpha^{(8)}_{2,6}   \left( \bm{p}_2 \otimes \bm{p}_6 \right)
            + \alpha^{(8)}_{4,4} \left( \bm{p}_4 \otimes \bm{p}_4 \right)
            + \alpha^{(8)}_{6,2} \left( \bm{p}_6 \otimes \bm{p}_2 \right)
        \text{.}
\end{equation}

When $h_* \to 0$,
the matrix $\bm{A}$ is decreasing faster than the rhs
$\bm{b}$, and therefore, the solution $\bm{w}$ must grow on the order of at least $h_*^{-2}$. The
sum of these increasing weights equals to one, and therefore, some of them must be negative. The
weights $\bm{w}$ may grow faster than $h_*^{-2}$ if such rapidly growing terms are orthogonal to
$\bm{b}$, to keep the $L^2$ finite. In Subsection~\labelcref{sec:SmallH:MM} we demonstrate that the
lowest power of $h_*$ is $-2M$. Thus, the weights are decomposed as:
\begin{equation}
    \label{eq:decw:h0:Odd}
        \bm{w} = \sum_{m = -M}^{+\infty} \bm{w}_{2 m} h_*^{2 m} \text{.}
\end{equation}

As the vector set $\bm{p}_{2m}$ does not form an orthonormal basis, we introduce complementary
vectors $\tilde{\bm{p}}_{2m}$, $m = 1, ..., M$, that are orthogonal to all $\bm{p}_{2j}$ for $j <
M$ and $j \ne m$. The complementary vectors are normalised such that $\langle \tilde{\bm{p}}_{2m},
\bm{p}_{2m} \rangle = 1$, where the angular brackets denote the inner product. Like $\bm{p}_{2m}$, the $\tilde{\bm{p}}_{2m}$ vectors are linearly
independent, but not mutually orthogonal. They are no longer orthogonal to $\bm{p}_{2j}$ when $j >
M$. For example, if $M = 3$, the space of weight vectors is $3$-D, the vectors $\bm{p}_{2k}$, $k = 1,
2, 3$, are not parallel, and vectors $\tilde{\bm{p}}_{2k}$, $k = 1, 2, 3$, always exist. However, a
vector $\tilde{\bm{p}}_8$ orthogonal to all $\bm{p}_k$, $k = 2, 4, 6$, cannot exist.

We represent the coefficients $\bm{w}_{2m}$ in the power series \labelcref{eq:decw:h0:Odd} of
mixand weights as linear combinations of the complementary vectors:
\begin{equation}
    \label{eq:decw2m:h0:Odd}
        \bm{w}_{2 m} = \sum_{j=1}^M C_{2j}^{(2m)} \tilde{\bm{p}}_{2j} \text{.}
\end{equation}
This form of decomposition is instrumental in combination with matrices expressed through outer products in
\cref{eq:A:h0:Odd} because for arbitrary vectors $\bm{x}$, $\bm{y}$, and $\bm{z}$ we have
$\left( \bm{x} \otimes \bm{y} \right) \bm{z} = \bm{x} \langle \bm{y}, \bm{z} \rangle$.
If $\bm{y}$ is one the vectors from
\cref{eq:p2n:Odd} and $\bm{z}$ is one of the complementary vectors $\tilde{\bm{p}}_{2m}$, we can use
their orthogonality to simplify the product of $\bm{A}$ and $\bm{w}$.

We substitute the decompositions
\labelcref{eq:deca:h0:Odd,eq:decb:h0:Odd,eq:decw:h0:Odd,eq:decw2m:h0:Odd} into the linear system
$\bm{A} \bm{w} = \bm{b}$, equate the terms with the same powers of $h_*$ on the lhs and the rhs, and then solve the
obtained linear systems of equations to find $C_{2j}^{(2m)}$. Using these results
in the product $\bm{w}{\tran} \bm{b}$, we then find the coefficients in the power series of $L^2$ (repeating
\cref{G1:eq:decL2:h0:Odd} here for completeness):
\begin{equation}
    \label{eq:decL2:h0:Odd}
    L^2 = \sum_{j=0}^{+\infty} L^2_{2j} h_*^{2j} \text{.}
\end{equation}

We first present
the solution up to the second-order terms in $L^2$ for the
simplest cases of $M = 1, 2$ and $3$, then generalise to arbitrary $M$ and arbitrary order of
decomposition in $h_*$.

\subsubsection{\texorpdfstring{$M = 1$}{M = 1}}\label{sec:SmallH:M1}

In this degenerate case $\bm{A}$ is not a matrix, but a scalar denoted as $A$. Likewise, $\bm{b}
\Rightarrow b$ and $\bm{w} \Rightarrow w$. All vectors $\bm{p}_{2j}$ become the scalar $1$, while
the ``orthogonal'' vectors $\tilde{\bm{p}}_{2j}$ cannot be defined, making the $M = 1$ case an
exception in terms of notation. Applying the series decomposition, the $h_*^2$ terms yield
\begin{equation}
    \label{eq:hstar2:m1}
    A_4 w_{-2} = \alpha^{(4)}_{2,2} w_{-2} = \beta_2
\end{equation}
to find $w_{-2}$. We skip the notation $C^{(2m)}_2$ because these values are equal to $w_{2m}$,
$\forall m \geq -1$. From the $h_*^4$ terms we obtain
\begin{equation}
    \label{eq:hstar4:m1}
    \alpha^{(4)}_{2,2} w_{0} + \left( \alpha^{(6)}_{2,4} + \alpha^{(6)}_{4,2} \right) w_{-2} = \beta_4
\end{equation}
and solve for $w_{0}$. Then the lowest coefficients in the series of $L^2$ are found as
\begin{equation}
    \label{eq:l2hstar0:m1}
    L^2_0 = {\lVert F \rVert}_2 - \beta_2 w_{-2}
\end{equation}
and
\begin{equation}
    \label{eq:l2hstar2:m1}
    L^2_2 = - \beta_2 w_{0} - \beta_4 w_{-2} \text{.}
\end{equation}
In \cref{eq:l2hstar0:m1}, ${\lVert F \rVert}_2$ is the $L^2$ norm squared of the standard Gaussian
function $\widetilde{\mathcal{N}} = \mathcal{N}(x; 0, 1)$ minus the central mixand $\mathcal{N}_0 = \mathcal{N}(x; 0,
\sigma)$, see \cref{G1:eq:cd:Odd}:
\begin{equation}
    \label{eq:F2}
        {\lVert F \rVert}_2
            = \int_{-\infty}^{+\infty} {\left( \widetilde{\mathcal{N}} - \mathcal{N}_0 \right)}^2 dx
            = \frac{1}{2 \sqrt{\pi}} - \frac{2}{\sqrt{2 \pi (1 + \sigma^2)}} + \frac{1}{2 \sigma \sqrt{\pi}} \text{.}
\end{equation}

\subsubsection{\texorpdfstring{$M = 2$}{M = 2}}\label{sec:SmallH:M2}

We define the system of equations from the series decomposition of $\bm{A} \bm{w} = \bm{b}$:
\begin{equation}
    \label{eq:Awb:h0:Odd:2}
    \begin{alignedat}{11}
    &\text{the $h_*^{0}$ terms} \quad \quad &&\bm{A}_4 \bm{w}_{-4} &&                       &&                       &&                          &&= \bm{0}   \\
    &\text{the $h_*^{2}$ terms} \quad \quad &&\bm{A}_4 \bm{w}_{-2} &&+ \bm{A}_6 \bm{w}_{-4} &&                       &&                          &&= \bm{b}_2 \\
    &\text{the $h_*^{4}$ terms} \quad \quad &&\bm{A}_4 \bm{w}_{0}  &&+ \bm{A}_6 \bm{w}_{-2} &&+ \bm{A}_8 \bm{w}_{-4} &&                          &&= \bm{b}_4 \\
    &\text{the $h_*^{6}$ terms} \quad \quad &&\bm{A}_4 \bm{w}_{2}  &&+ \bm{A}_6 \bm{w}_{0}  &&+ \bm{A}_8 \bm{w}_{-2} &&+ \bm{A}_{10} \bm{w}_{-4} &&= \bm{b}_6 \text{.}
    \end{alignedat}
\end{equation}
Expanding the $h_*^0$ equation gives
\begin{equation}
    \label{eq:hstar0:m2}
    \begin{aligned}[b]
        &\bm{A}_4 \bm{w}_{-4} = \alpha^{(4)}_{2, 2} (\bm{p}_2 \otimes \bm{p}_2) (C_2^{(-4)} \tilde{\bm{p}}_2 + C_4^{(-4)} \tilde{\bm{p}}_4) \\
        &= \alpha^{(4)}_{2, 2} (C_2^{(-4)} \langle \bm{p}_2, \tilde{\bm{p}}_2 \rangle \bm{p}_2 +
                                C_4^{(-4)} \langle \bm{p}_2, \tilde{\bm{p}}_4 \rangle \bm{p}_2) \\
        &= \alpha^{(4)}_{2, 2} (C_2^{(-4)} (1) \bm{p}_2 + C_4^{(-4)} (0) \bm{p}_2) \\
        &= \alpha^{(4)}_{2, 2} C_2^{(-4)} \bm{p}_2 = \bm{0} \text{.}
    \end{aligned}
\end{equation}
Hence, $C_2^{(-4)} = 0$, while $C_4^{(-4)}$ is undefined. Then, the $h_*^{2}$ equation is:
\begin{equation}
    \notag
    \begin{aligned}[b]
        &\bm{A}_4 \bm{w}_{-2} + \bm{A}_6 \bm{w}_{-4} \\
        &= \alpha^{(4)}_{2, 2} (C_2^{(-2)} \langle \bm{p}_2, \tilde{\bm{p}}_2 \rangle \bm{p}_2 +
                                C_4^{(-2)} \langle \bm{p}_2, \tilde{\bm{p}}_4 \rangle \bm{p}_2) + \\
        &\quad \alpha^{(6)}_{2, 4} C_4^{(-4)} \langle \bm{p}_4, \tilde{\bm{p}}_4 \rangle \bm{p}_2 +
               \alpha^{(6)}_{4, 2} C_4^{(-4)} \langle \bm{p}_2, \tilde{\bm{p}}_4 \rangle \bm{p}_4 \\
        &= \alpha^{(4)}_{2, 2} C_2^{(-2)} \bm{p}_2 +
           \alpha^{(6)}_{2, 4} C_4^{(-4)} \bm{p}_2 = \beta_2 \bm{p}_2 \text{,}
    \end{aligned}
\end{equation}
which yields
\begin{equation}
    \label{eq:hstar2:m2}
    \alpha^{(4)}_{2, 2} C_2^{(-2)} + \alpha^{(6)}_{2, 4} C_4^{(-4)} = \beta_2 \text{.}
\end{equation}
Now, moving on to the $h_*^{4}$ equation, and noting that for $M = 2$, the vectors
$\tilde{\bm{p}}_2$ and $\tilde{\bm{p}}_4$ are only orthogonal to $\bm{p}_4$ and $\bm{p}_2$,
respectively, not $\bm{p}_{2m}$ for $m \geq 3$, we obtain
\begin{align}
    \notag
    &\bm{A}_4 \bm{w}_{0} + \bm{A}_6 \bm{w}_{-2} + \bm{A}_8 \bm{w}_{-4} \\
    \notag
    \begin{split}
    \!\begin{alignedat}[b]{2}
        &= \ &&\alpha^{(4)}_{2, 2} (C_2^{(0)} \langle \bm{p}_2, \tilde{\bm{p}}_2 \rangle \bm{p}_2 +
                                    C_4^{(0)} \langle \bm{p}_2, \tilde{\bm{p}}_4 \rangle \bm{p}_2) + \\
        & &&\alpha^{(6)}_{2, 4} (C_2^{(-2)} \langle \bm{p}_4, \tilde{\bm{p}}_2 \rangle \bm{p}_2 +
                                 C_4^{(-2)} \langle \bm{p}_4, \tilde{\bm{p}}_4 \rangle \bm{p}_2) + \\
        & &&\alpha^{(6)}_{4, 2} (C_2^{(-2)} \langle \bm{p}_2, \tilde{\bm{p}}_2 \rangle \bm{p}_4 +
                                 C_4^{(-2)} \langle \bm{p}_2, \tilde{\bm{p}}_4 \rangle \bm{p}_4) + \\
        & &&\alpha^{(8)}_{2, 6} C_4^{(-4)} \langle \bm{p}_6, \tilde{\bm{p}}_4 \rangle \bm{p}_2 +
            \alpha^{(8)}_{4, 4} C_4^{(-4)} \langle \bm{p}_4, \tilde{\bm{p}}_4 \rangle \bm{p}_4 +
            \alpha^{(8)}_{6, 2} C_4^{(-4)} \langle \bm{p}_2, \tilde{\bm{p}}_4 \rangle \bm{p}_6 \\
    \end{alignedat}
    \end{split} \\
    \label{eq:hstar4:m2}
    \begin{split}
    \!\begin{alignedat}[b]{2}
        &= \ &&\alpha^{(4)}_{2, 2} C_2^{(0)} \bm{p}_2 +
               \alpha^{(6)}_{2, 4} C_4^{(-2)} \bm{p}_2 +
               \alpha^{(6)}_{4, 2} C_2^{(-2)} \bm{p}_4 + \\
        & &&\alpha^{(8)}_{2, 6} C_4^{(-4)} \langle \bm{p}_6, \tilde{\bm{p}}_4 \rangle \bm{p}_2 +
             \alpha^{(8)}_{4, 4} C_4^{(-4)} \bm{p}_4 = \beta_4 \bm{p}_4 \text{.}
    \end{alignedat}
    \end{split}
\end{align}
Multiplying both sides by $\tilde{\bm{p}}_4$ yields
\begin{equation}
    \label{eq:hstar4:p4:m2}
    \alpha^{(6)}_{4, 2} C_2^{(-2)} + \alpha^{(8)}_{4, 4} C_4^{(-4)} = \beta_4 \text{.}
\end{equation}
Combined with \cref{eq:hstar2:m2}, this gives a linear system for $C_2^{(-2)}$ and $C_4^{(-4)}$.
Solving it, we find the term of order $h_*^0$ in the power series decomposition of $L^2$:
\begin{equation}
    \label{eq:l2hstar0:m2}
    \begin{aligned}[b]
        L^2_0 &= {\lVert F \rVert}_2 - \bm{w}_{-2}\tran \bm{b}_{2} - \bm{w}_{-4}\tran \bm{b}_{4} \\
              &= {\lVert F \rVert}_2 - \beta_2 (C_2^{(-2)} \langle \bm{p}_2, \tilde{\bm{p}}_2 \rangle +
                                                C_4^{(-2)} \langle \bm{p}_2, \tilde{\bm{p}}_4 \rangle)
                                     - \beta_4 C_4^{(-4)} \langle \bm{p}_4, \tilde{\bm{p}}_4 \rangle \\
              &= {\lVert F \rVert}_2 - \beta_2 C_2^{(-2)} - \beta_4 C_4^{(-4)} \text{,}
    \end{aligned}
\end{equation}
where ${\lVert F \rVert}_2$ is given in \cref{eq:F2}.
To get the coefficients for the second-order term $L^2_2$, we first multiply \cref{eq:hstar4:m2} by
$\tilde{\bm{p}}_2$ to obtain
\begin{equation}
    \label{eq:hstar4:p2:m2}
    \alpha^{(4)}_{2, 2} C_2^{(0)} + \alpha^{(6)}_{2, 4} C_4^{(-2)}
    = -\alpha^{(8)}_{2, 6} C_4^{(-4)} \langle \bm{p}_6, \tilde{\bm{p}}_4 \rangle \text{,}
\end{equation}
Then we consider the $h_*^6$ equation:
\begin{align}
    \notag
    &\bm{A}_4 \bm{w}_{2} + \bm{A}_6 \bm{w}_{0} + \bm{A}_8 \bm{w}_{-2} + \bm{A}_{10} \bm{w}_{-4} \\
    \notag
    \begin{split}
    \!\begin{alignedat}[b]{2}
        &= \ &&\alpha^{(4)}_{2, 2} (C_2^{(2)} \langle \bm{p}_2, \tilde{\bm{p}}_2 \rangle \bm{p}_2 +
                                    C_4^{(2)} \langle \bm{p}_2, \tilde{\bm{p}}_4 \rangle \bm{p}_2) + \\
        & &&\alpha^{(6)}_{2, 4} (C_2^{(0)} \langle \bm{p}_4, \tilde{\bm{p}}_2 \rangle \bm{p}_2 +
                                 C_4^{(0)} \langle \bm{p}_4, \tilde{\bm{p}}_4 \rangle \bm{p}_2) + \\
        & &&\alpha^{(6)}_{4, 2} (C_2^{(0)} \langle \bm{p}_2, \tilde{\bm{p}}_2 \rangle \bm{p}_4 +
                                 C_4^{(0)} \langle \bm{p}_2, \tilde{\bm{p}}_4 \rangle \bm{p}_4) + \\
        & &&\alpha^{(8)}_{2, 6} (C_2^{(-2)} \langle \bm{p}_6, \tilde{\bm{p}}_2 \rangle \bm{p}_2 +
                                 C_4^{(-2)} \langle \bm{p}_6, \tilde{\bm{p}}_4 \rangle \bm{p}_2) + \\
        & &&\alpha^{(8)}_{4, 4} (C_2^{(-2)} \langle \bm{p}_4, \tilde{\bm{p}}_2 \rangle \bm{p}_4 +
                                 C_4^{(-2)} \langle \bm{p}_4, \tilde{\bm{p}}_4 \rangle \bm{p}_4) + \\
        & &&\alpha^{(8)}_{6, 2} (C_2^{(-2)} \langle \bm{p}_2, \tilde{\bm{p}}_2 \rangle \bm{p}_6 +
                                 C_4^{(-2)} \langle \bm{p}_2, \tilde{\bm{p}}_4 \rangle \bm{p}_6) + \\
        & &&\alpha^{(10)}_{2, 8} C_4^{(-4)} \langle \bm{p}_8, \tilde{\bm{p}}_4 \rangle \bm{p}_2 +
            \alpha^{(10)}_{4, 6} C_4^{(-4)} \langle \bm{p}_6, \tilde{\bm{p}}_4 \rangle \bm{p}_4 + \\
        & &&\alpha^{(10)}_{6, 4} C_4^{(-4)} \langle \bm{p}_4, \tilde{\bm{p}}_4 \rangle \bm{p}_6 +
            \alpha^{(10)}_{8, 2} C_4^{(-4)} \langle \bm{p}_2, \tilde{\bm{p}}_4 \rangle \bm{p}_8 \\
    \end{alignedat}
    \end{split} \\
    \label{eq:hstar6:m2}
    \begin{split}
    \!\begin{alignedat}[b]{2}
        &= \ &&\alpha^{(4)}_{2, 2} C_2^{(2)} \bm{p}_2 +
               \alpha^{(6)}_{2, 4} C_4^{(0)} \bm{p}_2 +
               \alpha^{(6)}_{4, 2} C_2^{(0)} \bm{p}_4 + \\
        & &&\alpha^{(8)}_{2, 6} (C_2^{(-2)} \langle \bm{p}_6, \tilde{\bm{p}}_2 \rangle \bm{p}_2 +
                                 C_4^{(-2)} \langle \bm{p}_6, \tilde{\bm{p}}_4 \rangle \bm{p}_2) \\
        & &&\alpha^{(8)}_{4, 4} C_4^{(-2)} \bm{p}_4 +
            \alpha^{(8)}_{6, 2} C_2^{(-2)} \bm{p}_6 +
            \alpha^{(10)}_{2, 8} C_4^{(-4)} \langle \bm{p}_8, \tilde{\bm{p}}_4 \rangle \bm{p}_2 + \\
        & &&\alpha^{(10)}_{4, 6} C_4^{(-4)} \langle \bm{p}_6, \tilde{\bm{p}}_4 \rangle \bm{p}_4 +
            \alpha^{(10)}_{6, 4} C_4^{(-4)}  \bm{p}_6 = \beta_6 \bm{p}_6 \text{.}
    \end{alignedat}
    \end{split}
\end{align}
Multiplying by $\tilde{\bm{p}}_4$ and combining the common terms, we simplify the above to
\begin{align}
    \label{eq:hstar6:p4:m2}
    \begin{split}
    \!\begin{alignedat}[b]{2}
    &\alpha^{(6)}_{4, 2} C_2^{(0)} + \alpha^{(8)}_{4, 4} C_4^{(-2)} \\
    &= \left( \beta_{6} - \alpha^{(10)}_{4, 6} C_4^{(-4)} - \alpha^{(8)}_{6, 2} C^{(-2)}_2 - \alpha^{(10)}_{6, 4} C_4^{(-4)} \right) \langle \bm{p}_6, \tilde{\bm{p}}_4 \rangle \text{,}
    \end{alignedat}
    \end{split}
\end{align}
Together with \cref{eq:hstar4:p2:m2}, we solve for $C_2^{(0)}$ and $C_4^{(-2)}$. Finally,
we take the $h_*^2$ terms in $\bm{w}{\tran} \bm{b}$ to obtain
\begin{equation}
    \label{eq:l2hstar2:m2}
    \begin{aligned}[b]
        L^2_2 &= -\bm{w}_{0}\tran \bm{b}_2 - \bm{w}_{-2}\tran \bm{b}_4 - \bm{w}_{-4}\tran \bm{b}_6 \\
              &= -\beta_2 (C_2^{(0)} \langle \bm{p}_2, \tilde{\bm{p}}_2 \rangle +
                           C_4^{(0)} \langle \bm{p}_2, \tilde{\bm{p}}_4 \rangle) \\
              &\quad -\beta_4 (C_2^{(-2)} \langle \bm{p}_4, \tilde{\bm{p}}_2 \rangle +
                               C_4^{(-2)} \langle \bm{p}_4, \tilde{\bm{p}}_4 \rangle) \\
              &\quad -\beta_6 C_4^{(-4)} \langle \bm{p}_6, \tilde{\bm{p}}_4 \rangle \\
              &= -\beta_2 C_2^{(0)} - \beta_4 C_4^{(-2)}
                 -\beta_6 C_4^{(-4)} \langle \bm{p}_6, \tilde{\bm{p}}_4 \rangle \text{.}
    \end{aligned}
\end{equation}
The order of solution for various coefficients $C_{2j}^{(2m)}$ is given in \cref{tbl:m2}.

\begin{table}[ht]
    \caption{Order of appearance of various unknowns in the equations of the system
    \labelcref{eq:Awb:h0:Odd:2}. Bright pink box shows the $h_*^0$-order equation for
    $C^{(-4)}_2$. Bright green boxes indicate the system for the coefficients used in $L^2_0$, and bright
    yellow boxes show the corresponding system for the coefficients used in $L^2_2$. Pale colours mark the
    variables already known from the elements in the corresponding bright colour. The rhs
    (last column) differentiates only between zero and non-zero (check mark) entries.}
    \label{tbl:m2}
    \begin{center}
    \begin{tabular}{ l l c c c c c c c c }
        order   & vector     & $C^{(-4)}_2$ & $C^{(-4)}_4$ & $C^{(-2)}_2$ & $C^{(-2)}_4$ & $C^{(0)}_2$ & $C^{(0)}_4$ & $C^{(2)}_2$ & rhs        \\
        \hline
        \hline
        $h_*^0$ & $\bm{p}_2$ & \crc{35}     & -            & -            & -            & -           & -           & -           & \crz{35}   \\
        \hline
        $h_*^2$ & $\bm{p}_2$ & -            & \cgc{45}     & \cgc{45}     & -            & -           & -           & -           & \cgc{45}   \\
                & $\bm{p}_4$ & \crc{7}      & -            & -            & -            & -           & -           & -           & $\bm{0}$   \\
        \hline
        $h_*^4$ & $\bm{p}_2$ & \crc{7}      & \cgc{10}     & -            & \cyc{55}     & \cyc{55}    & -           & -           & \cyz{55}   \\
                & $\bm{p}_4$ & -            & \cgc{45}     & \cgc{45}     & -            & -           & -           & -           & \cgc{45}   \\
                & $\bm{p}_6$ & \crc{7}      & -            & -            & -            & -           & -           & -           & $\bm{0}$   \\
        \hline
        $h_*^6$ & $\bm{p}_2$ & \crc{7}      & \cgc{10}     & \cgc{10}     & \cyc{15}     & -           & \Checkmark  & \Checkmark  & $\bm{0}$   \\
                & $\bm{p}_4$ & \crc{7}      & \cgc{10}     & \cgc{10}     & \cyc{55}     & \cyc{55}    & -           & -           & \cyz{55}   \\
                & $\bm{p}_6$ & -            & \cgc{10}     & \cgc{10}     & -            & -           & -           & -           & \Checkmark \\
                & $\bm{p}_8$ & \crc{7}      & -            & -            & -            & -           & -           & -           & $\bm{0}$   \\
        \hline
        \hline
    \end{tabular}
    \end{center}
\end{table}

\subsubsection{\texorpdfstring{$M = 3$}{M = 3}}\label{sec:SmallH:M3}

Taking the equal powers of $h_*$ in the decomposition of $\bm{A} \bm{w} = \bm{b}$, we obtain
\begin{equation}
\begin{alignedat}[c]{13}
    \label{eq:Awb:h0:Odd:3}
        &\bm{A}_4 \bm{w}_{-6} &&                       &&                       &&                          &&                          &&  &&= \bm{0} \\
        &\bm{A}_4 \bm{w}_{-4} &&+ \bm{A}_6 \bm{w}_{-6} &&                       &&                          &&                          &&  &&= \bm{0} \\
        &\bm{A}_4 \bm{w}_{-2} &&+ \bm{A}_6 \bm{w}_{-4} &&+ \bm{A}_8 \bm{w}_{-6} &&                          &&                          &&  &&= \bm{b}_2 \\
        &\bm{A}_4 \bm{w}_{0}  &&+ \bm{A}_6 \bm{w}_{-2} &&+ \bm{A}_8 \bm{w}_{-4} &&+ \bm{A}_{10} \bm{w}_{-6} &&                          &&  &&= \bm{b}_4 \\
        &\bm{A}_4 \bm{w}_{2}  &&+ \bm{A}_6 \bm{w}_{0}  &&+ \bm{A}_8 \bm{w}_{-2} &&+ \bm{A}_{10} \bm{w}_{-4} &&+ \bm{A}_{12} \bm{w}_{-6} &&  &&= \bm{b}_6 \\
        &\bm{A}_4 \bm{w}_{4}  &&+ \bm{A}_6 \bm{w}_{2}  &&+ \bm{A}_8 \bm{w}_{0}  &&+ \bm{A}_{10} \bm{w}_{-2} &&+ \bm{A}_{12} \bm{w}_{-4} &&+ \bm{A}_{14} \bm{w}_{-6} &&= \bm{b}_8 \text{,}
\end{alignedat}
\end{equation}
where the top equation corresponds to $h_*^{-2}$, the second to $h_*^0$, and then $h_*^{2m}$ for
every $\bm{b}_{2m}$ in the rhs.

Starting with the $h_*^{-2}$ equation, we find
\begin{equation}
    \label{eq:hstar-2:m3}
    \begin{aligned}[b]
    &\bm{A}_{4} \bm{w}_{-6}
    = \alpha^{(4)}_{2, 2} (C_2^{(-6)} \langle \bm{p}_2, \tilde{\bm{p}}_2 \rangle \bm{p}_2 +
                           C_4^{(-6)} \langle \bm{p}_2, \tilde{\bm{p}}_4 \rangle \bm{p}_2 +
                           C_6^{(-6)} \langle \bm{p}_2, \tilde{\bm{p}}_6 \rangle \bm{p}_2) \\
    &= \alpha^{(4)}_{2, 2} C_2^{(-6)} \bm{p}_2 = \bm{0} \text{.}
    \end{aligned}
\end{equation}
Hence, $C_2^{(-6)} = 0$, while $C_4^{(-6)}$ and $C_6^{(-6)}$ are undefined. Then, the $h_*^{0}$
equation yields
\begin{equation}
    \notag
    \begin{aligned}[b]
    &\bm{A}_{4} \bm{w}_{-4} + \bm{A}_{6} \bm{w}_{-6} \\
    &= \alpha^{(4)}_{2, 2} (C_2^{(-4)} \langle \bm{p}_2, \tilde{\bm{p}}_2 \rangle \bm{p}_2 +
                            C_4^{(-4)} \langle \bm{p}_2, \tilde{\bm{p}}_4 \rangle \bm{p}_2 +
                            C_6^{(-4)} \langle \bm{p}_2, \tilde{\bm{p}}_6 \rangle \bm{p}_2) + \\
    &\quad \alpha^{(6)}_{2, 4} (C_4^{(-6)} \langle \bm{p}_4, \tilde{\bm{p}}_4 \rangle \bm{p}_2 +
                                C_6^{(-6)} \langle \bm{p}_4, \tilde{\bm{p}}_6 \rangle \bm{p}_2) + \\
    &\quad \alpha^{(6)}_{4, 2} (C_4^{(-6)} \langle \bm{p}_2, \tilde{\bm{p}}_4 \rangle \bm{p}_4 +
                                C_6^{(-6)} \langle \bm{p}_2, \tilde{\bm{p}}_6 \rangle \bm{p}_4) \\
    &= \alpha^{(4)}_{2, 2} C_2^{(-4)} \bm{p}_2 +
       \alpha^{(6)}_{2, 4} C_4^{(-6)} \bm{p}_2 = \bm{0} \text{,}
    \end{aligned}
\end{equation}
or simply
\begin{equation}
    \label{eq:hstar0:m3}
    \alpha^{(4)}_{2, 2} C_2^{(-4)} + \alpha^{(6)}_{2, 4} C_4^{(-6)} = 0 \text{.}
\end{equation}
Moving on to the $h_*^{2}$ equation, we obtain
\begin{align}
    \notag
    \bm{A}_{4} &\bm{w}_{-2} + \bm{A}_{6} \bm{w}_{-4} + \bm{A}_{8} \bm{w}_{-6} \\
    \notag
    \begin{split}
    \!\begin{alignedat}[b]{2}
        &= \ &&\alpha^{(4)}_{2, 2} (C_2^{(-2)} \langle \bm{p}_2, \tilde{\bm{p}}_2 \rangle \bm{p}_2 +
                                C_4^{(-2)} \langle \bm{p}_2, \tilde{\bm{p}}_4 \rangle \bm{p}_2 +
                                C_6^{(-2)} \langle \bm{p}_2, \tilde{\bm{p}}_6 \rangle \bm{p}_2) + \\
        & &&\alpha^{(6)}_{2, 4} (C_2^{(-4)} \langle \bm{p}_4, \tilde{\bm{p}}_2 \rangle \bm{p}_2 +
                                  C_4^{(-4)} \langle \bm{p}_4, \tilde{\bm{p}}_4 \rangle \bm{p}_2 +
                                  C_6^{(-4)} \langle \bm{p}_4, \tilde{\bm{p}}_6 \rangle \bm{p}_2) + \\
        & &&\alpha^{(6)}_{4, 2} (C_2^{(-4)} \langle \bm{p}_2, \tilde{\bm{p}}_2 \rangle \bm{p}_4 +
                                  C_4^{(-4)} \langle \bm{p}_2, \tilde{\bm{p}}_4 \rangle \bm{p}_4 +
                                  C_6^{(-4)} \langle \bm{p}_2, \tilde{\bm{p}}_6 \rangle \bm{p}_4) + \\
        & &&\alpha^{(8)}_{2, 6} (C_4^{(-6)} \langle \bm{p}_6, \tilde{\bm{p}}_4 \rangle \bm{p}_2 +
                                  C_6^{(-6)} \langle \bm{p}_6, \tilde{\bm{p}}_6 \rangle \bm{p}_2) + \\
        & &&\alpha^{(8)}_{4, 4} (C_4^{(-6)} \langle \bm{p}_4, \tilde{\bm{p}}_4 \rangle \bm{p}_4 +
                                  C_6^{(-6)} \langle \bm{p}_4, \tilde{\bm{p}}_6 \rangle \bm{p}_4) + \\
        & &&\alpha^{(8)}_{6, 2} (C_4^{(-6)} \langle \bm{p}_2, \tilde{\bm{p}}_4 \rangle \bm{p}_6 +
                                  C_6^{(-6)} \langle \bm{p}_2, \tilde{\bm{p}}_6 \rangle \bm{p}_6) \\
    \end{alignedat}
    \end{split}
    \\
    \label{eq:hstar2:m3}
    \begin{split}
    \!\begin{alignedat}[b]{2}
        &= \ &&\alpha^{(4)}_{2, 2} C_2^{(-2)} \bm{p}_2 +
           \alpha^{(6)}_{2, 4} C_4^{(-4)} \bm{p}_2 +
           \alpha^{(6)}_{4, 2} C_2^{(-4)} \bm{p}_4 + \\
        & &&\alpha^{(8)}_{2, 6} C_6^{(-6)} \bm{p}_2 +
             \alpha^{(8)}_{4, 4} C_4^{(-6)} \bm{p}_4 = \beta_2 \bm{p}_2 \text{.}
    \end{alignedat}
    \end{split}
\end{align}
Multiplied by $\tilde{\bm{p}}_4$, \cref{eq:hstar2:m3} reduces to
\begin{equation}
    \label{eq:hstar2:p4:m3}
    \alpha^{(6)}_{4, 2} C_2^{(-4)} + \alpha^{(8)}_{4, 4} C_4^{(-6)} = \bm{0} \text{,}
\end{equation}
which, combined with \cref{eq:hstar0:m3}, yields a system for $C_2^{(-4)}$ and $C_4^{(-6)}$, with
the solution $C_2^{(-4)} = C_4^{(-6)} = 0$. If instead we multiply \cref{eq:hstar2:m3} by
$\tilde{\bm{p}}_2$, it yields
\begin{equation}
    \label{eq:hstar2:p2:m3}
    \alpha^{(4)}_{2, 2} C_2^{(-2)} +
    \alpha^{(6)}_{2, 4} C_4^{(-4)} +
    \alpha^{(8)}_{2, 6} C_6^{(-6)}  = \beta_2 \text{.}
\end{equation}
Now consider the $h_*^{4}$ equation:
\begin{align}
    \notag
    &\bm{A}_{4} \bm{w}_{0} + \bm{A}_{6} \bm{w}_{-2} + \bm{A}_{8} \bm{w}_{-4} + \bm{A}_{10} \bm{w}_{-6} \\
    \notag
    \begin{split}
    \!\begin{alignedat}[b]{2}
    &= \ &&\alpha^{(4)}_{2, 2} (C_2^{(0)} \langle \bm{p}_2, \tilde{\bm{p}}_2 \rangle \bm{p}_2 +
                                C_4^{(0)} \langle \bm{p}_2, \tilde{\bm{p}}_4 \rangle \bm{p}_2 +
                                C_6^{(0)} \langle \bm{p}_2, \tilde{\bm{p}}_6 \rangle \bm{p}_2) + \\
    & &&\alpha^{(6)}_{2, 4} (C_2^{(-2)} \langle \bm{p}_4, \tilde{\bm{p}}_2 \rangle \bm{p}_2 +
                             C_4^{(-2)} \langle \bm{p}_4, \tilde{\bm{p}}_4 \rangle \bm{p}_2 +
                             C_6^{(-2)} \langle \bm{p}_4, \tilde{\bm{p}}_6 \rangle \bm{p}_2) + \\
    & &&\alpha^{(6)}_{4, 2} (C_2^{(-2)} \langle \bm{p}_2, \tilde{\bm{p}}_2 \rangle \bm{p}_4 +
                             C_4^{(-2)} \langle \bm{p}_2, \tilde{\bm{p}}_4 \rangle \bm{p}_4 +
                             C_6^{(-2)} \langle \bm{p}_2, \tilde{\bm{p}}_6 \rangle \bm{p}_4) + \\
    & &&\alpha^{(8)}_{2, 6} (C_4^{(-4)} \langle \bm{p}_6, \tilde{\bm{p}}_4 \rangle \bm{p}_2 +
                             C_6^{(-4)} \langle \bm{p}_6, \tilde{\bm{p}}_6 \rangle \bm{p}_2) + \\
    & &&\alpha^{(8)}_{4, 4} (C_4^{(-4)} \langle \bm{p}_4, \tilde{\bm{p}}_4 \rangle \bm{p}_4 +
                             C_6^{(-4)} \langle \bm{p}_4, \tilde{\bm{p}}_6 \rangle \bm{p}_4) + \\
    & &&\alpha^{(8)}_{6, 2} (C_4^{(-4)} \langle \bm{p}_2, \tilde{\bm{p}}_4 \rangle \bm{p}_6 +
                             C_6^{(-4)} \langle \bm{p}_2, \tilde{\bm{p}}_6 \rangle \bm{p}_6) + \\
    & &&\alpha^{(10)}_{2, 8} C_6^{(-6)} \langle \bm{p}_8, \tilde{\bm{p}}_6 \rangle \bm{p}_2 +
        \alpha^{(10)}_{4, 6} C_6^{(-6)} \langle \bm{p}_6, \tilde{\bm{p}}_6 \rangle \bm{p}_4 + \\
    & &&\alpha^{(10)}_{6, 4} C_6^{(-6)} \langle \bm{p}_4, \tilde{\bm{p}}_6 \rangle \bm{p}_6 +
        \alpha^{(10)}_{8, 2} C_6^{(-6)} \langle \bm{p}_2, \tilde{\bm{p}}_6 \rangle \bm{p}_8 \\
    \end{alignedat}
    \end{split} \\
    \label{eq:hstar4:m3}
    \begin{split}
    \!\begin{alignedat}[b]{2}
    &= \ &&\alpha^{(4)}_{2, 2} C_2^{(0)} \bm{p}_2 +
           \alpha^{(6)}_{2, 4} C_4^{(-2)} \bm{p}_2 +
           \alpha^{(6)}_{4, 2} C_2^{(-2)} \bm{p}_4 + \\
    & &&\alpha^{(8)}_{2, 6} C_6^{(-4)} \bm{p}_2 +
        \alpha^{(8)}_{4, 4} C_4^{(-4)} \bm{p}_4 +
        \alpha^{(10)}_{2, 8} C_6^{(-6)} \langle \bm{p}_8, \tilde{\bm{p}}_6 \rangle \bm{p}_2 + \\
    & &&\alpha^{(10)}_{4, 6} C_6^{(-6)} \bm{p}_4 = \beta_4 \bm{p}_4 \text{.}
    \end{alignedat}
    \end{split}
\end{align}
Multiplication by $\tilde{\bm{p}}_4$ reduces \cref{eq:hstar4:m3} to:
\begin{equation}
    \label{eq:hstar4:p4:m3}
         \alpha^{(6)}_{4, 2} C_2^{(-2)}  +
         \alpha^{(8)}_{4, 4} C_4^{(-4)}  +
         \alpha^{(10)}_{4, 6} C_6^{(-6)} = \beta_4 \text{.}
\end{equation}
Going next to the $h_*^{6}$ equation, we find
\begin{align}
    \notag
    &\bm{A}_{4} \bm{w}_{2} + \bm{A}_{6} \bm{w}_{0} + \bm{A}_{8} \bm{w}_{-2} + \bm{A}_{10} \bm{w}_{-4} + \bm{A}_{12} \bm{w}_{-6} \\
    \notag
    \begin{split}
    \!\begin{alignedat}[b]{2}
    &= \ &&\alpha^{(4)}_{2, 2} (C_2^{(2)} \langle \bm{p}_2, \tilde{\bm{p}}_2 \rangle \bm{p}_2 +
                                C_4^{(2)} \langle \bm{p}_2, \tilde{\bm{p}}_4 \rangle \bm{p}_2 +
                                C_6^{(2)} \langle \bm{p}_2, \tilde{\bm{p}}_6 \rangle \bm{p}_2) + \\
    & &&\alpha^{(6)}_{2, 4} (C_2^{(0)} \langle \bm{p}_4, \tilde{\bm{p}}_2 \rangle \bm{p}_2 +
                             C_4^{(0)} \langle \bm{p}_4, \tilde{\bm{p}}_4 \rangle \bm{p}_2 +
                             C_6^{(0)} \langle \bm{p}_4, \tilde{\bm{p}}_6 \rangle \bm{p}_2) + \\
    & &&\alpha^{(6)}_{4, 2} (C_2^{(0)} \langle \bm{p}_2, \tilde{\bm{p}}_2 \rangle \bm{p}_4 +
                             C_4^{(0)} \langle \bm{p}_2, \tilde{\bm{p}}_4 \rangle \bm{p}_4 +
                             C_6^{(0)} \langle \bm{p}_2, \tilde{\bm{p}}_6 \rangle \bm{p}_4) + \\
    & &&\alpha^{(8)}_{2, 6} (C_2^{(-2)} \langle \bm{p}_6, \tilde{\bm{p}}_2 \rangle \bm{p}_2 +
                             C_4^{(-2)} \langle \bm{p}_6, \tilde{\bm{p}}_4 \rangle \bm{p}_2 +
                             C_6^{(-2)} \langle \bm{p}_6, \tilde{\bm{p}}_6 \rangle \bm{p}_2) + \\
    & &&\alpha^{(8)}_{4, 4} (C_2^{(-2)} \langle \bm{p}_4, \tilde{\bm{p}}_2 \rangle \bm{p}_4 +
                             C_4^{(-2)} \langle \bm{p}_4, \tilde{\bm{p}}_4 \rangle \bm{p}_4 +
                             C_6^{(-2)} \langle \bm{p}_4, \tilde{\bm{p}}_6 \rangle \bm{p}_4) + \\
    & &&\alpha^{(8)}_{6, 2} (C_2^{(-2)} \langle \bm{p}_2, \tilde{\bm{p}}_2 \rangle \bm{p}_6 +
                             C_4^{(-2)} \langle \bm{p}_2, \tilde{\bm{p}}_4 \rangle \bm{p}_6 +
                             C_6^{(-2)} \langle \bm{p}_2, \tilde{\bm{p}}_6 \rangle \bm{p}_6) + \\
    & &&\alpha^{(10)}_{2, 8} (C_4^{(-4)} \langle \bm{p}_8, \tilde{\bm{p}}_4 \rangle \bm{p}_2 +
                              C_6^{(-4)} \langle \bm{p}_8, \tilde{\bm{p}}_6 \rangle \bm{p}_2) + \\
    & &&\alpha^{(10)}_{4, 6} (C_4^{(-4)} \langle \bm{p}_6, \tilde{\bm{p}}_4 \rangle \bm{p}_4 +
                              C_6^{(-4)} \langle \bm{p}_6, \tilde{\bm{p}}_6 \rangle \bm{p}_4) + \\
    & &&\alpha^{(10)}_{6, 4} (C_4^{(-4)} \langle \bm{p}_4, \tilde{\bm{p}}_4 \rangle \bm{p}_6 +
                              C_6^{(-4)} \langle \bm{p}_4, \tilde{\bm{p}}_6 \rangle \bm{p}_6) + \\
    & &&\alpha^{(10)}_{8, 2} (C_4^{(-4)} \langle \bm{p}_2, \tilde{\bm{p}}_4 \rangle \bm{p}_8 +
                              C_6^{(-4)} \langle \bm{p}_2, \tilde{\bm{p}}_6 \rangle \bm{p}_8) + \\
    & &&\alpha^{(12)}_{2, 10} C_6^{(-6)} \langle \bm{p}_{10}, \tilde{\bm{p}}_6 \rangle \bm{p}_2 +
        \alpha^{(12)}_{4, 8} C_6^{(-6)} \langle \bm{p}_{8}, \tilde{\bm{p}}_6 \rangle \bm{p}_4 + \\
    & &&\alpha^{(12)}_{6, 6} C_6^{(-6)} \langle \bm{p}_{6}, \tilde{\bm{p}}_6 \rangle \bm{p}_6 +
        \alpha^{(12)}_{8, 4} C_6^{(-6)} \langle \bm{p}_{4}, \tilde{\bm{p}}_6 \rangle \bm{p}_8 + \\
    & &&\alpha^{(12)}_{10, 2} C_6^{(-6)} \langle \bm{p}_{2}, \tilde{\bm{p}}_6 \rangle \bm{p}_{10}
    \end{alignedat}
    \end{split}
    \end{align}
    \vspace{-0.5cm}
    \begin{align}
    \label{eq:hstar6:m3}
    \begin{split}
    \!\begin{alignedat}[b]{2}
    &= \ &&\alpha^{(4)}_{2, 2} C_2^{(2)} \bm{p}_2 +
           \alpha^{(6)}_{2, 4} C_4^{(0)} \bm{p}_2 +
           \alpha^{(6)}_{4, 2} C_2^{(0)} \bm{p}_4 + \\
    & &&\alpha^{(8)}_{2, 6} C_6^{(-2)} \bm{p}_2 +
        \alpha^{(8)}_{4, 4} C_4^{(-2)} \bm{p}_4 +
        \alpha^{(8)}_{6, 2} C_2^{(-2)} \bm{p}_6 + \\
    & &&\alpha^{(10)}_{2, 8} (C_4^{(-4)} \langle \bm{p}_8, \tilde{\bm{p}}_4 \rangle \bm{p}_2 +
                              C_6^{(-4)} \langle \bm{p}_8, \tilde{\bm{p}}_6 \rangle \bm{p}_2) + \\
    & &&\alpha^{(10)}_{4, 6} C_6^{(-4)} \bm{p}_4 + \alpha^{(10)}_{6, 4} C_4^{(-4)} \bm{p}_6 +
        \alpha^{(12)}_{2, 10} C_6^{(-6)} \langle \bm{p}_{10}, \tilde{\bm{p}}_6 \rangle \bm{p}_2 + \\
    & &&\alpha^{(12)}_{4, 8} C_6^{(-6)} \langle \bm{p}_{8}, \tilde{\bm{p}}_6 \rangle \bm{p}_4 +
        \alpha^{(12)}_{6, 6} C_6^{(-6)} \bm{p}_6 = \beta_6 \bm{p}_6 \text{.}
    \end{alignedat}
    \end{split}
\end{align}
Multiplying by $\tilde{\bm{p}}_6$, we obtain
\begin{equation}
    \label{eq:hstar6:p6:m3}
         \alpha^{(8)}_{6, 2} C_2^{(-2)} +
         \alpha^{(10)}_{6, 4} C_4^{(-4)} +
         \alpha^{(12)}_{6, 6} C_6^{(-6)} = \beta_4 \text{.}
\end{equation}
Together with \cref{eq:hstar2:p2:m3,eq:hstar4:p4:m3}, this gives a system of equations for $C_2^{(-2)}$, $C_4^{(-4)}$,
and $C_6^{(-6)}$, and allows to find $L^2_0$:
\begin{equation}
    \label{eq:l2hstar0:m3}
    \begin{aligned}[b]
        L^2_0 &= {\lVert F \rVert}_2 - \bm{w}_{-2}\tran \bm{b}_{2} - \bm{w}_{-4}\tran \bm{b}_{4} - \bm{w}_{-6}\tran \bm{b}_{6} \\
              &= {\lVert F \rVert}_2 - \beta_2 (C_2^{(-2)} \langle \bm{p}_2, \tilde{\bm{p}}_2 \rangle +
                                                C_4^{(-2)} \langle \bm{p}_2, \tilde{\bm{p}}_4 \rangle +
                                                C_6^{(-2)} \langle \bm{p}_2, \tilde{\bm{p}}_6 \rangle) \\
              &\quad                 - \beta_4 (C_4^{(-4)} \langle \bm{p}_4, \tilde{\bm{p}}_4 \rangle +
                                                C_6^{(-4)} \langle \bm{p}_4, \tilde{\bm{p}}_6 \rangle) \\
              &\quad                 - \beta_6 C_6^{(-6)} \langle \bm{p}_6, \tilde{\bm{p}}_6 \rangle \\
              &= {\lVert F \rVert}_2 - \beta_2 C_2^{(-2)} - \beta_4 C_4^{(-4)} - \beta_6 C_6^{(-6)} \text{,}
    \end{aligned}
\end{equation}
where ${\lVert F \rVert}_2$ is given in \cref{eq:F2}.

To find the next batch of coefficients, we multiply \cref{eq:hstar4:m3} by $\tilde{\bm{p}}_2$,
yielding
\begin{equation}
    \label{eq:hstar4:p2:m3}
    \alpha^{(4)}_{2, 2} C_2^{(0)} + \alpha^{(6)}_{2, 4} C_4^{(-2)} + \alpha^{(8)}_{2, 6} C_6^{(-4)}
    = -\alpha^{(10)}_{2, 8} C_6^{(-6)} \langle \bm{p}_8, \tilde{\bm{p}}_6 \rangle \text{.}
\end{equation}
Similarly, multiplication of \cref{eq:hstar6:m3} by $\tilde{\bm{p}}_4$ gives:
\begin{equation}
    \label{eq:hstar6:p4:m3}
    \alpha^{(6)}_{4, 2} C_2^{(0)} + \alpha^{(8)}_{4, 4} C_4^{(-2)} + \alpha^{(10)}_{4, 6} C_6^{(-4)}
    = -\alpha^{(12)}_{4, 8} C_6^{(-6)} \langle \bm{p}_8, \tilde{\bm{p}}_6 \rangle \text{.}
\end{equation}
Finally, we consider the $h_*^8$ equation:
\begin{align}
    \notag
    &\bm{A}_{4} \bm{w}_{4} + \bm{A}_{6} \bm{w}_{2} + \bm{A}_{8} \bm{w}_{0} + \bm{A}_{10} \bm{w}_{-2} + \bm{A}_{12} \bm{w}_{-4} +
     \bm{A}_{12} \bm{w}_{-6} \\
    \notag
    \begin{split}
    \!\begin{alignedat}[b]{2}
    &= \ &&\alpha^{(4)}_{2, 2} (C_2^{(4)} \langle \bm{p}_2, \tilde{\bm{p}}_2 \rangle \bm{p}_2 +
                                C_4^{(4)} \langle \bm{p}_2, \tilde{\bm{p}}_4 \rangle \bm{p}_2 +
                                C_6^{(4)} \langle \bm{p}_2, \tilde{\bm{p}}_6 \rangle \bm{p}_2) + \\
    & &&\alpha^{(6)}_{2, 4} (C_2^{(2)} \langle \bm{p}_4, \tilde{\bm{p}}_2 \rangle \bm{p}_2 +
                             C_4^{(2)} \langle \bm{p}_4, \tilde{\bm{p}}_4 \rangle \bm{p}_2 +
                             C_6^{(2)} \langle \bm{p}_4, \tilde{\bm{p}}_6 \rangle \bm{p}_2) + \\
    & &&\alpha^{(6)}_{4, 2} (C_2^{(2)} \langle \bm{p}_2, \tilde{\bm{p}}_2 \rangle \bm{p}_4 +
                             C_4^{(2)} \langle \bm{p}_2, \tilde{\bm{p}}_4 \rangle \bm{p}_4 +
                             C_6^{(2)} \langle \bm{p}_2, \tilde{\bm{p}}_6 \rangle \bm{p}_4) + \\
    & &&\alpha^{(8)}_{2, 6} (C_2^{(0)} \langle \bm{p}_6, \tilde{\bm{p}}_2 \rangle \bm{p}_2 +
                             C_4^{(0)} \langle \bm{p}_6, \tilde{\bm{p}}_4 \rangle \bm{p}_2 +
                             C_6^{(0)} \langle \bm{p}_6, \tilde{\bm{p}}_6 \rangle \bm{p}_2) + \\
    & &&\alpha^{(8)}_{4, 4} (C_2^{(0)} \langle \bm{p}_4, \tilde{\bm{p}}_2 \rangle \bm{p}_4 +
                             C_4^{(0)} \langle \bm{p}_4, \tilde{\bm{p}}_4 \rangle \bm{p}_4 +
                             C_6^{(0)} \langle \bm{p}_4, \tilde{\bm{p}}_6 \rangle \bm{p}_4) + \\
    & &&\alpha^{(8)}_{6, 2} (C_2^{(0)} \langle \bm{p}_2, \tilde{\bm{p}}_2 \rangle \bm{p}_6 +
                             C_4^{(0)} \langle \bm{p}_2, \tilde{\bm{p}}_4 \rangle \bm{p}_6 +
                             C_6^{(0)} \langle \bm{p}_2, \tilde{\bm{p}}_6 \rangle \bm{p}_6) + \\
    & &&\alpha^{(10)}_{2, 8} (C_2^{(-2)} \langle \bm{p}_8, \tilde{\bm{p}}_2 \rangle \bm{p}_2 +
                              C_4^{(-2)} \langle \bm{p}_8, \tilde{\bm{p}}_4 \rangle \bm{p}_2 +
                              C_6^{(-2)} \langle \bm{p}_8, \tilde{\bm{p}}_6 \rangle \bm{p}_2) + \\
    & &&\alpha^{(10)}_{4, 6} (C_2^{(-2)} \langle \bm{p}_6, \tilde{\bm{p}}_2 \rangle \bm{p}_4 +
                              C_4^{(-2)} \langle \bm{p}_6, \tilde{\bm{p}}_4 \rangle \bm{p}_4 +
                              C_6^{(-2)} \langle \bm{p}_6, \tilde{\bm{p}}_6 \rangle \bm{p}_4) + \\
    & &&\alpha^{(10)}_{6, 4} (C_2^{(-2)} \langle \bm{p}_4, \tilde{\bm{p}}_2 \rangle \bm{p}_6 +
                              C_4^{(-2)} \langle \bm{p}_4, \tilde{\bm{p}}_4 \rangle \bm{p}_6 +
                              C_6^{(-2)} \langle \bm{p}_4, \tilde{\bm{p}}_6 \rangle \bm{p}_6) + \\
    & &&\alpha^{(10)}_{8, 2} (C_2^{(-2)} \langle \bm{p}_2, \tilde{\bm{p}}_2 \rangle \bm{p}_8 +
                              C_4^{(-2)} \langle \bm{p}_2, \tilde{\bm{p}}_4 \rangle \bm{p}_8 +
                              C_6^{(-2)} \langle \bm{p}_2, \tilde{\bm{p}}_6 \rangle \bm{p}_8) + \\
    & &&\alpha^{(12)}_{2, 10} (C_4^{(-4)} \langle \bm{p}_{10}, \tilde{\bm{p}}_4 \rangle \bm{p}_2 +
                               C_6^{(-4)} \langle \bm{p}_{10}, \tilde{\bm{p}}_6 \rangle \bm{p}_2) + \\
    & &&\alpha^{(12)}_{4, 8} (C_4^{(-4)} \langle \bm{p}_8, \tilde{\bm{p}}_4 \rangle \bm{p}_4 +
                              C_6^{(-4)} \langle \bm{p}_8, \tilde{\bm{p}}_6 \rangle \bm{p}_4) + \\
    & &&\alpha^{(12)}_{6, 6} (C_4^{(-4)} \langle \bm{p}_6, \tilde{\bm{p}}_4 \rangle \bm{p}_6 +
                              C_6^{(-4)} \langle \bm{p}_6, \tilde{\bm{p}}_6 \rangle \bm{p}_6) + \\
    & &&\alpha^{(12)}_{8, 4} (C_4^{(-4)} \langle \bm{p}_4, \tilde{\bm{p}}_4 \rangle \bm{p}_8 +
                              C_6^{(-4)} \langle \bm{p}_4, \tilde{\bm{p}}_6 \rangle \bm{p}_8) + \\
    & &&\alpha^{(12)}_{10, 2} (C_4^{(-4)} \langle \bm{p}_2, \tilde{\bm{p}}_4 \rangle \bm{p}_{10} +
                               C_6^{(-4)} \langle \bm{p}_2, \tilde{\bm{p}}_6 \rangle \bm{p}_{10}) + \\
    & &&\alpha^{(14)}_{2, 12} C_6^{(-6)} \langle \bm{p}_{12}, \tilde{\bm{p}}_6 \rangle \bm{p}_2 +
        \alpha^{(14)}_{4, 10} C_6^{(-6)} \langle \bm{p}_{10}, \tilde{\bm{p}}_6 \rangle \bm{p}_4 + \\
    & &&\alpha^{(14)}_{6, 8} C_6^{(-6)} \langle \bm{p}_8, \tilde{\bm{p}}_6 \rangle \bm{p}_6 +
        \alpha^{(14)}_{8, 6} C_6^{(-6)} \langle \bm{p}_6, \tilde{\bm{p}}_6 \rangle \bm{p}_8 + \\
    & &&\alpha^{(14)}_{10, 4} C_6^{(-6)} \langle \bm{p}_4, \tilde{\bm{p}}_6 \rangle \bm{p}_{10} +
        \alpha^{(14)}_{12, 2} C_6^{(-6)} \langle \bm{p}_2, \tilde{\bm{p}}_6 \rangle \bm{p}_{12} \\
    \end{alignedat}
    \end{split} \\
    \label{eq:hstar8:m3}
    \begin{split}
    \begin{alignedat}[b]{2}
    &= \ &&\alpha^{(4)}_{2, 2} C_2^{(4)} \bm{p}_2 +
           \alpha^{(6)}_{2, 4} C_4^{(2)} \bm{p}_2 +
           \alpha^{(6)}_{4, 2} C_2^{(2)} \bm{p}_4 + \\
    & &&\alpha^{(8)}_{2, 6} C_6^{(0)} \bm{p}_2 +
        \alpha^{(8)}_{4, 4} C_4^{(0)} \bm{p}_4 +
        \alpha^{(8)}_{6, 2} C_2^{(0)} \bm{p}_6 + \\
    & &&\alpha^{(10)}_{2, 8} (C_2^{(-2)} \langle \bm{p}_8, \tilde{\bm{p}}_2 \rangle \bm{p}_2 +
                              C_4^{(-2)} \langle \bm{p}_8, \tilde{\bm{p}}_4 \rangle \bm{p}_2 +
                              C_6^{(-2)} \langle \bm{p}_8, \tilde{\bm{p}}_6 \rangle \bm{p}_2) + \\
    & &&\alpha^{(10)}_{4, 6} C_6^{(-2)} \bm{p}_4 +
        \alpha^{(10)}_{6, 4} C_4^{(-2)} \bm{p}_6 +
        \alpha^{(10)}_{8, 2} C_2^{(-2)} \bm{p}_8 + \\
    & &&\alpha^{(12)}_{2, 10} (C_4^{(-4)} \langle \bm{p}_{10}, \tilde{\bm{p}}_4 \rangle \bm{p}_2 +
                               C_6^{(-4)} \langle \bm{p}_{10}, \tilde{\bm{p}}_6 \rangle \bm{p}_2) + \\
    & &&\alpha^{(12)}_{4, 8} (C_4^{(-4)} \langle \bm{p}_8, \tilde{\bm{p}}_4 \rangle \bm{p}_4 +
                              C_6^{(-4)} \langle \bm{p}_8, \tilde{\bm{p}}_6 \rangle \bm{p}_4) + \\
    & &&\alpha^{(12)}_{6, 6} C_6^{(-4)} \bm{p}_6 +
        \alpha^{(12)}_{8, 4} C_4^{(-4)} \bm{p}_8 +
        \alpha^{(14)}_{2, 12} C_6^{(-6)} \langle \bm{p}_{12}, \tilde{\bm{p}}_6 \rangle \bm{p}_2 + \\
    & &&\alpha^{(14)}_{4, 10} C_6^{(-6)} \langle \bm{p}_{10}, \tilde{\bm{p}}_6 \rangle \bm{p}_4 +
        \alpha^{(14)}_{6, 8} C_6^{(-6)} \langle \bm{p}_8, \tilde{\bm{p}}_6 \rangle \bm{p}_6 + \\
    & &&\alpha^{(14)}_{8, 6} C_6^{(-6)} \bm{p}_8 = \beta_8 \bm{p}_8 \text{.}
    \end{alignedat}
    \end{split}
\end{align}
Multiplying \cref{eq:hstar8:m3} by $\tilde{\bm{p}}_6$ and putting known terms to the right, we get:
\begin{align}
    \label{eq:hstar8:p6:m3}
    \begin{split}
    \!\begin{alignedat}[b]{2}
    &\alpha^{(8)}_{6, 2} C_2^{(0)} + \alpha^{(10)}_{6, 4} C_4^{(-2)} + \alpha^{(12)}_{6, 6} C_6^{(-4)} \\
    &= \left( \beta_{8} - \alpha^{(14)}_{6, 8} C_6^{(-6)}  - \alpha^{(10)}_{8, 2} C_2^{(-2)} -
       \alpha^{(12)}_{8, 4} C_4^{(-4)} - \alpha^{(14)}_{8, 6} C_6^{(-6)}\right) \langle \bm{p}_8, \tilde{\bm{p}}_6 \rangle \text{.}
    \end{alignedat}
    \end{split}
\end{align}
Combined with \cref{eq:hstar4:p2:m3,eq:hstar6:p4:m3}, we solve for $C_2^{(0)}$, $C_4^{(-2)}$,
and $C_6^{(-4)}$, which allows us to obtain $L^2_2$:
\begin{equation}
    \label{eq:l2hstar2:m3}
    \begin{aligned}[b]
        L^2_2 &= -\bm{w}_{0}\tran \bm{b}_2 - \bm{w}_{-2}\tran \bm{b}_4 - \bm{w}_{-4}\tran \bm{b}_6 - \bm{w}_{-6}\tran \bm{b}_8 \\
              &= -\beta_2 (C_2^{(0)} \langle \bm{p}_2, \tilde{\bm{p}}_2 \rangle +
                           C_4^{(0)} \langle \bm{p}_2, \tilde{\bm{p}}_4 \rangle +
                           C_6^{(0)} \langle \bm{p}_2, \tilde{\bm{p}}_6 \rangle) \\
                 &\quad -\beta_4 (C_2^{(-2)} \langle \bm{p}_4, \tilde{\bm{p}}_2 \rangle +
                           C_4^{(-2)} \langle \bm{p}_4, \tilde{\bm{p}}_4 \rangle +
                           C_6^{(-2)} \langle \bm{p}_4, \tilde{\bm{p}}_6 \rangle) \\
                 &\quad -\beta_6 (C_4^{(-4)} \langle \bm{p}_6, \tilde{\bm{p}}_4 \rangle +
                           C_6^{(-4)} \langle \bm{p}_6, \tilde{\bm{p}}_6 \rangle) \\
                 &\quad -\beta_8 C_6^{(-6)} \langle \bm{p}_8, \tilde{\bm{p}}_6 \rangle \\
              &= -\beta_2 C_2^{(0)} - \beta_4 C_4^{(-2)}
                 -\beta_6 C_6^{(-4)} - \beta_8 C_6^{(-6)} \langle \bm{p}_8, \tilde{\bm{p}}_6 \rangle \text{.}
    \end{aligned}
\end{equation}
The order of solution for various coefficients $C_{2j}^{(2m)}$ is given in \cref{tbl:m3}. Together
with \cref{tbl:m2}, they establish a pattern how the equation systems for isolated subsets of
$C_{2j}^{(2m)}$ appear in consecutive groups spanning several orders of the $h_*$ powers. This forms
the basis of the general $M$ solution presented next.

\begin{sidewaystable}[ht]
    \caption{Order of appearance of various unknowns in the equations of the system
    \labelcref{eq:Awb:h0:Odd:3}. Bright pink box shows the $h_*^0$-order equation for
    $C^{(-6)}_2$. Blue boxes show the system of \cref{eq:hstar0:m3} and \cref{eq:hstar2:p4:m3}.
    Green and yellow boxes indicate the systems for the coefficients used in $L^2_0$ and $L^2_2$,
    respectively.
    The values of $C^{(2)}_4$, $C^{(2)}_6$,
    and $C^{(2)}_2$ that appear only in \cref{eq:hstar8:m3}, the highest-order equation considered,
    were omitted for clarity. Pale colours mark the variables already known from the elements
    in the corresponding bright colour. The rhs (last column) differentiates only
    between zero and non-zero (check mark) entries.}
    \label{tbl:m3}
\begin{center}
\begin{tabular}{ l l c c c c c c c c c c c c c c }
    order      & vector        & $C^{(-6)}_2$ & $C^{(-6)}_4$ & $C^{(-6)}_6$ & $C^{(-4)}_2$ & $C^{(-4)}_4$ & $C^{(-4)}_6$ & $C^{(-2)}_2$ & $C^{(-2)}_4$ & $C^{(-2)}_6$ & $C^{(0)}_2$ & $C^{(0)}_4$ & $C^{(0)}_6$ & $C^{(2)}_2$ & rhs        \\
    \hline
    \hline
    $h_*^{-2}$ & $\bm{p}_2$    & \crc{35}     & -            & -            & -            & -            & -            & -            & -            & -            & -           & -           & -           & -           & \crz{35}   \\
    \hline
    $h_*^0$    & $\bm{p}_2$    & -            & \cbc{35}     & -            & \cbc{35}     & -            & -            & -            & -            & -            & -           & -           & -           & -           & \cbz{35}   \\
               & $\bm{p}_4$    & \crc{7}      & -            & -            & -            & -            & -            & -            & -            & -            & -           & -           & -           & -           & $\bm{0}$   \\
    \hline
    $h_*^2$    & $\bm{p}_2$    & -            & -            & \cgc{45}     & -            & \cgc{45}     & -            & \cgc{45}     & -            & -            & -           & -           & -           & -           & \cgc{45}   \\
               & $\bm{p}_4$    & -            & \cbc{35}     & -            & \cbc{35}     & -            & -            & -            & -            & -            & -           & -           & -           & -           & \cbz{35}   \\
               & $\bm{p}_6$    & \crc{7}      & -            & -            & -            & -            & -            & -            & -            & -            & -           & -           & -           & -           & $\bm{0}$   \\
    \hline
    $h_*^4$    & $\bm{p}_2$    & \crc{7}      & \cbc{10}     & \cgc{10}     & -            & -            & \cyc{55}     & -            & \cyc{55}     & -            & \cyc{55}    & -           & -           & -           & \cyz{55}   \\
               & $\bm{p}_4$    & -            & -            & \cgc{45}     & -            & \cgc{45}     & -            & \cgc{45}     & -            & -            & -           & -           & -           & -           & \cgc{45}   \\
               & $\bm{p}_6$    & -            & \cbc{10}     & -            & -            & -            & -            & -            & -            & -            & -           & -           & -           & -           & $\bm{0}$   \\
               & $\bm{p}_8$    & \crc{7}      & -            & -            & \cbc{10}     & -            & -            & -            & -            & -            & -           & -           & -           & -           & $\bm{0}$   \\
    \hline
    $h_*^6$    & $\bm{p}_2$    & \crc{7}      & \cbc{10}     & \cgc{10}     & \cbc{10}     & \cgc{10}     & \cyc{15}     & -            & -            & \Checkmark   & -           & \Checkmark  & -           & \Checkmark  & $\bm{0}$   \\
               & $\bm{p}_4$    & \crc{7}      & \cbc{10}     & \cgc{10}     & -            & -            & \cyc{55}     & -            & \cyc{55}     & -            & \cyc{55}    & -           & -           & -           & \cyz{55}   \\
               & $\bm{p}_6$    & -            & -            & \cgc{45}     & -            & \cgc{45}     & -            & \cgc{45}     & -            & -            & -           & -           & -           & -           & \cgc{45}   \\
               & $\bm{p}_8$    & -            & \cbc{10}     & -            & \cbc{10}     & -            & -            & -            & -            & -            & -           & -           & -           & -           & $\bm{0}$   \\
               & $\bm{p}_{10}$ & \crc{7}      & -            & -            & -            & -            & -            & -            & -            & -            & -           & -           & -           & -           & $\bm{0}$   \\
    \hline
    $h_*^8$    & $\bm{p}_2$    & \crc{7}      & \cbc{10}     & \cgc{10}     & \cbc{10}     & \cgc{10}     & \cyc{15}     & \cgc{10}     & \cyc{15}     & \Checkmark   & -           & -           & \Checkmark  & -           & $\bm{0}$   \\
               & $\bm{p}_4$    & \crc{7}      & \cbc{10}     & \cgc{10}     & \cbc{10}     & \cgc{10}     & \cyc{15}     & -            & -            & \Checkmark   & -           & \Checkmark  & -           & \Checkmark  & $\bm{0}$   \\
               & $\bm{p}_6$    & \crc{7}      & \cbc{10}     & \cgc{10}     & -            & \cgc{10}     & \cyc{55}     & \cgc{10}     & \cyc{55}     & -            & \cyc{55}    & -           & -           & -           & \cyz{55}   \\
               & $\bm{p}_8$    & -            & -            & \cgc{10}     & -            & \cgc{10}     & -            & \cgc{10}     & -            & -            & -           & -           & -           & -           & \Checkmark \\
               & $\bm{p}_{10}$ & -            & \cbc{10}     & -            & \cbc{10}     & -            & -            & -            & -            & -            & -           & -           & -           & -           & $\bm{0}$   \\
               & $\bm{p}_{12}$ & \crc{7}      & -            & -            & -            & -            & -            & -            & -            & -            & -           & -           & -           & -           & $\bm{0}$   \\
    \hline
    \hline
\end{tabular}
\end{center}
\end{sidewaystable}

\subsubsection{The solution for general \texorpdfstring{$M$}{M}}\label{sec:SmallH:GenM}

Having analysed the special cases of $M = 1, 2, 3$, we now seek a solution for arbitrary dimensionality $M$.
Furthermore, we seek to prove the statement made before \cref{eq:decw:h0:Odd} that the highest power
of the decomposition of the weights in \cref{eq:decw:h0:Odd} is $-2M$.

\paragraph{Solution for the weights}

Using decompositions \labelcref{eq:decb:h0:Odd,eq:deca:h0:Odd,eq:decw:h0:Odd} and matching the multipliers of
$h_*^{2q}$ on both sides of $\bm{A} \bm{b} = \bm{w}$ for $q = -M + 2, -M + 3, \dots$, we obtain
\begin{equation}
    \label{eq:genM:Awbh:Odd}
    \sum_{j=2}^{q + M} \bm{A}_{2j} \bm{w}_{2(q - j)} = \beta_{2 q} \bm{p}_{2 q} [ q \geq 1 ] \text{,}
\end{equation}
where $[\ldots]$ denotes the Iverson bracket \cite[p.~24]{Graham:2017:Concrete}; it equals $1$ if $q
\ge 1$ and zero otherwise.
The rhs of \cref{eq:genM:2q:2n:Odd} uses non-existing
$\beta_{2q}$ factors for $q < 1$; however, for such cases the Iverson bracket evaluates to zero, so,
we ignore this inconsistency.

After expanding $\bm{A}_{2j}$ and $\bm{w}_{2(q - j)}$ into their own summations by \cref{eq:A:h0:Odd,eq:decw:h0:Odd},
the lhs of \cref{eq:genM:Awbh:Odd} becomes
\begin{equation}
    \label{eq:genM:2q:Odd}
    \begin{aligned}[b]
    \sum_{j=2}^{q + M} \bm{A}_{2j} \bm{w}_{2(q - j)}
    &= \sum_{j=2}^{q + M}
         \sum_{k=1}^{j-1} \alpha^{(2j)}_{2k, 2(j - k)} (\bm{p}_{2k} \otimes \bm{p}_{2(j - k)})
            \sum_{l=1}^M C^{(2(q - j))}_{2l} \tilde{\bm{p}}_{2l} \\
    &= \sum_{j=2}^{q + M} \sum_{k=1}^{j-1} \sum_{l=1}^M
         \alpha^{(2j)}_{2k, 2(j - k)} C^{(2(q - j))}_{2l}
         \langle \bm{p}_{2(j - k)}, \tilde{\bm{p}}_{2l} \rangle \bm{p}_{2k} \\
    &= \sum_{j=2}^{q + M} \sum_{i=1}^{j-1} \sum_{l=1}^M
         \alpha^{(2j)}_{2(j - i), 2i} C^{(2(q - j))}_{2l}
         \langle \bm{p}_{2i}, \tilde{\bm{p}}_{2l} \rangle \bm{p}_{2(j - i)} \text{,}
    \end{aligned}
\end{equation}
where on the last line we replaced the middle sum index as $k = j - i$.

Multiplying both sides of \cref{eq:genM:Awbh:Odd} by some $\tilde{\bm{p}}_{2n}$, we
finally get
\begin{multline}
    \label{eq:genM:2q:2n:Odd}
        \sum_{j=2}^{q + M} \sum_{i=1}^{j-1} \sum_{l=1}^M
            \alpha^{(2j)}_{2(j - i), 2i} C^{(2(q - j))}_{2l}
            \langle \bm{p}_{2i}, \tilde{\bm{p}}_{2l} \rangle
            \langle \bm{p}_{2(j - i)}, \tilde{\bm{p}}_{2n} \rangle \\
        = \beta_{2q} \langle \bm{p}_{2q}, \tilde{\bm{p}}_{2n} \rangle [q \geq 1]
        \text{.}
\end{multline}
We refer to \cref{eq:genM:2q:2n:Odd} for some $q, n$ as ``the equation $2q
/ 2n$'': it is the $h_*^{2q}$-th term in the power series of $\bm{A} \bm{b} = \bm{w}$ projected onto
$\tilde{\bm{p}}_{2n}$.

Consider the groups of $2q / 2n$ equations with $q = -M + s + n$ for $s \geq 1$ and $n$ running from
$1$ to $\min{(s,M)}$. The rhs of \cref{eq:genM:2q:2n:Odd} is non-zero when either $q = n
\le M$, which means $s = M$ (for all $n$), or $q > M$, which means $s + n > 2M$; the latter is
possible only when $s > M$. We will see that the $s = M$ group yields the equations for the
coefficients $C^{(2m)}_{2j}$ that are used in $L^2_0$. Similarly, solving the $s = M + 1$ group
produces the coefficients used in $L^2_2$. But before we get there, consider the groups for smaller
$s$.

\subparagraph{Group $s = 1$}

The simplest group for $s = 1$ has a single equation with $n = 1$, $q = -M + 2$ (the lowest value
possible), and, based on the sum limits in \cref{eq:genM:2q:2n:Odd}, $j = 2$ and $i = 1$. The first
inner product in \cref{eq:genM:2q:2n:Odd} becomes the Kronecker delta function $\delta_{1,l}$, and
the entire triple sum reduces to a single term:
\begin{equation}
    \label{eq:Seq1}
    \alpha^{(4)}_{2, 2} C^{(-2 M)}_{2} = \beta_{2(-M + 2)} \delta_{-M + 2,1} [-M + 2 \geq 1] \text{.}
\end{equation}
The rhs is non-zero only when $M = 1$, in which case we reproduce \cref{eq:hstar2:m1}.
Otherwise, if $M > 1$, we get $C^{(-2 M)}_{2} = 0$, also familiar from \cref{eq:hstar0:m2} for $M =
2$ and \cref{eq:hstar-2:m3} for $M = 3$.

\subparagraph{Groups $s < M$}

Now consider more general $s < M$.
By definitions of
$\bm{p}_{2k}$ and $\tilde{\bm{p}}_{2k}$, the inner products in \cref{eq:genM:2q:2n:Odd} satisfy the following conditions:
\begin{subequations}
    \label{eq:PInner}
\begin{align}
    \langle \bm{p}_{2i}, \tilde{\bm{p}}_{2l} \rangle &= \delta_{i,l} && \text{ if } 1 \leq i \leq M
    \label{eq:PInner:a} \text{,} \\
    \langle \bm{p}_{2i}, \tilde{\bm{p}}_{2l} \rangle &> 1            && \text{ if } i > M
    \label{eq:PInner:b} \text{.}
\end{align}
\end{subequations}
The second index is
$1 \le l \le M$.
We aim to use the Kronecker delta functions from \cref{eq:PInner:a} to
simplify the summation in \cref{eq:genM:2q:2n:Odd}. However, according to \cref{eq:PInner:b}, the
deltas might not be available. If we ignore this complication for a moment (it is reviewed below),
then the only non-zero terms in \cref{eq:genM:2q:2n:Odd} are those with $j = n + i$ and $l = i$. As
$i \ge 1$, we must have $j \ge n + 1$, increasing the lower limit of the first sum. Also, by
definition of $s$ and $n$, we have $q + M = n + s$. Together this yields
\begin{equation}
    \notag
    \sum_{j=n + 1}^{n + s}
        \biggl. \alpha^{(2j)}_{2(j - i), 2i} C^{(2(q - j))}_{2l} \biggr|_{
            \raisemath{3pt}{
                \begin{subarray}{l}
                l = i \\i = j - n
                \end{subarray}
            }
        } = 0
    \text{,}
\end{equation}
or, using a new summation index $m = j - n$,
\begin{equation}
    \label{eq:SlessM}
    \sum_{m=1}^{s}
        \alpha^{(2(m + n))}_{2n, 2m} C^{(2(-M + s - m))}_{2m} = 0 \quad
    \text{ for $n = 1, \ldots, s$.}
\end{equation}
Note that the indices of the $C^{(\ldots)}_{2m}$ variables in \cref{eq:SlessM} depend only on $m$,
$1 \leq m \leq s$, but not on $n$, so we have $s$ linear equations for the same $s$ unknowns. The
coefficients $\alpha^{(2(m + n))}_{2n, 2m}$ of this linear system are given in \cref{eq:alpha:h0:Odd}
and do not exhibit any linear dependencies across the rows or columns. Thus, the system is solvable,
and given that all the right-hand sides are zeros, so are the solutions.

For $M = 3$, the system \labelcref{eq:SlessM} at $s = 2$ yields \cref{eq:hstar0:m3,eq:hstar2:p4:m3},
while for $s = 1$ and $M \ne 1$ it becomes \cref{eq:Seq1}, which was considered separately because
of the special case of non-zero rhs at $M = 1$.

So far, however, the derivation tacitly omitted a possible hurdle: the indices of the $\bm{p}_{2i}$
and $\bm{p}_{2(j - i)}$ vectors in \cref{eq:genM:2q:2n:Odd} may exceed $M$, and if so, the compaction
of the two inner sums by Kronecker deltas would not apply. Now we examine if such indices may occur.

If we have large $n$, $j$ could be large enough so that the index $i$ in the first inner product in
\cref{eq:genM:2q:2n:Odd} exceeds $M$. Then the product would fall into the category of
\cref{eq:PInner:b}, not \cref{eq:PInner:a}. However, when $i > M$, we would have $j - i < j -
M \leq s + n - M < n$, as we consider $1 \le s < M$. Therefore, the second inner product would be
zero, so the overall contribution of such terms to the sum would be nil. Note that for $s = M$ we
would have $s + n - M = n$, but still $j - i < n$. However, for $s > M$, this condition would no
longer hold.

Consider now the second inner product. May we have $j - i > M$? Indeed, $\max{(j - i)} = s + n - 1$,
and it may exceed $M$ for sufficiently large $s > (M + 1) / 2$, as $\max{(n)} = s$.\footnote{Given
that $s$ is an integer and $s < M$, the condition may be satisfied only when $M \ge 4$.} The total
contribution of such terms to the lhs of \cref{eq:genM:2q:2n:Odd} would be
\begin{equation}
    \label{eq:NonKron}
    \sum_{j=M + 2}^{s + n} \sum_{i=1}^{j - M - 1}
        \alpha^{(2j)}_{2(j - i), 2i} C^{(2(q - j))}_{2i}
        \langle \bm{p}_{2(j - i)}, \tilde{\bm{p}}_{2n} \rangle
        \text{.}
\end{equation}
We adjusted the sum limits because to have $j - i > M$ (and $i \ge 1$), one has to have $j \ge M +
2$ and $i \leq j - M - 1$. Also, we already know that we may not have $i > M$ and $j - i > M$ at the
same time. Therefore, the first inner product is $\delta_{i,l}$. Expressing $q$ through $s$ and $n$,
\cref{eq:NonKron} contains the variables
\begin{equation}
    \label{eq:NonKron:C}
    C^{(2(-M + s + n - j))}_{2i} \quad \text{ for $M + 2 \le j \le s + n$, $1 \le i \le j - M - 1$.}
\end{equation}
If we solve the length-$s$ groups of equations sequentially from smaller to larger $s$, we would
already know all these values by solving \cref{eq:SlessM} for some $s' < s$. More formally, for all
$i, j, s, n$ such that
\begin{equation}
    \notag
    \begin{aligned}
        M + 2 & \le j \le s + n     \text{,}        \quad && 1 \le s \le M - 1 \text{,} \\
            1 & \le i \le j - M - 1 \text{,}        \quad && 1 \le n \le s     \text{,} \\
    \end{aligned}
\end{equation}
we can find $s'$ and $m = i \le s'$ that satisfy
\begin{equation}
    \notag
    1 \le s' < s
    \quad \text{and} \quad
    -M + s + n - j = -M + s' - m \text{.}
\end{equation}
The proof is straightforward. We explicitly have
\begin{equation}
    \notag
    s' = s + n - j + i \text{.}
\end{equation}
Then
\begin{subequations}
    \label{eq:Sdash:Lim}
\begin{equation}
    \label{eq:Sdash:Min}
    \min\limits_{i, j, n}{(s + n - j + i)} =
        \min\limits_{j, n}{(s + n - j + 1)} = 1 \text{,}
\end{equation}
as $\max{(j)} = s + n$, while
\begin{multline}
    \label{eq:Sdash:Max}
    \max\limits_{i, j, n}{(s + n - j + i)} =
        \max\limits_{i, j, n}{(s + n - (j - i))} = \\
        \max\limits_{n}{(s + n - (M + 1))} =
        2 s - (M + 1) < s \text{,}
\end{multline}
as $s < M$. Therefore, $1 \le s' < s$. At the same time, $m = i = s' - (s + n - j) \le s'$, because
$s + n \ge j$. Thus, all variables in \cref{eq:NonKron:C} would already appear in the linear systems
\labelcref{eq:SlessM} for $s = s'$, and consequently, would be zeros. This completes the
proof of \cref{eq:SlessM} validity: when some of the inner products in \cref{eq:genM:2q:2n:Odd} do
not evaluate to Kronecker delta functions, the corresponding terms are equal to zero for other reasons.

Also note that the conclusions given by \cref{eq:Sdash:Lim} hold for arbitrary $s \ge M$. The minimal
value in \cref{eq:Sdash:Min} remains the same, while for the maximal value we obtain a slightly
modified expression
\begin{equation}
    \label{eq:Sdash:Max:2}
    \max\limits_{i, j, n}{(s + n - j + i)} =
        \max\limits_{n}{(s + n - (M + 1))} =
        s - 1 < s \text{,}
\end{equation}
\end{subequations}
because the maximal value of $n$ would now be $M$.

\subparagraph{Group $s = M$}

This happens to be a minor extension of the earlier derivation for $s < M$. The main difference is
that for $s = M$ we have $q = n$, and therefore, both the inner product and the Iverson bracket in
the rhs of the $2q/2n$ equation \labelcref{eq:genM:2q:2n:Odd} are equal to one. The
presented earlier analysis of the inner products in the lhs remains the same: the
derivation applies verbatim for $s = M$. Thus, the $s = M$ group yields $M$ linear equations
\begin{equation}
    \label{eq:SeqM}
    \sum_{m=1}^{M}
        \alpha^{(2(m + n))}_{2n, 2m} C^{(-2m)}_{2m} = \beta_{2n} \quad
    \text{ for $n = 1, \ldots, M$.}
\end{equation}
This linear system generalises the earlier obtained \cref{eq:hstar2:m2,eq:hstar4:p4:m2} for $M =
2$, and \cref{eq:hstar2:p2:m3,eq:hstar4:p4:m3,eq:hstar6:p6:m3} for $M = 3$.

\subparagraph{Group $s = M + 1$}

Now both inner products in \cref{eq:genM:2q:2n:Odd} may fall into the category of \cref{eq:PInner:b}.
Thus, we need to consider three types of terms in the lhs of \cref{eq:genM:2q:2n:Odd}:
\begin{itemize}
    \item The ``regular'' terms with $i \le M$ and $j - i \le M$; for these, the inner products
        would evaluate to Kronecker deltas.
    \item The terms with $i > M$ and $j - i \le M$, so that the second inner product evaluates to
        $\delta_{j - i, n}$ for some $n$.
    \item The terms with $j - i > M$ and arbitrary $i$.
\end{itemize}
Also, the rhs of \cref{eq:genM:2q:2n:Odd} could be non-zero.

The first type of terms is simple and would eventually yield the lhs of a linear system
similar to \cref{eq:SlessM}. It is written out below.

For the second type, we have $j = n + i > M + n$. Given that the upper limit of $j$ is $M + n + 1$,
this leaves only $j = M + n + 1$, and therefore, $i = M + 1$, $j - i = n$, and $q - j = -M$. This
eliminates the first two summations leaving
\begin{equation}
    \notag
    \sum_{l = 1}^{M}
        \alpha^{(2(M + 1 + n))}_{2n, 2(M + 1)} C^{(-2M)}_{2l}
            \langle \bm{p}_{2(M + 1)}, \tilde{\bm{p}}_{2l} \rangle \text{.}
\end{equation}
As we know that $C^{(-2M)}_{2l} = 0$ unless $l = M$, the sum reduces to a single term
\begin{equation}
    \notag
    \alpha^{(2(M + 1 + n))}_{2n, 2(M + 1)} C^{(-2M)}_{2M}
        \langle \bm{p}_{2(M + 1)}, \tilde{\bm{p}}_{2M} \rangle \text{.}
\end{equation}
This component exists for all $n$.

The contribution of the third type of terms is given by \cref{eq:NonKron}. It was proven above
that all the variables $C^{(\ldots)}_{2i}$ in that expression are already known from solving
\cref{eq:SlessM} for some $s' < s$. Furthermore, only the variables $C^{(-2m)}_{2m}$ with matching
lower and upper indices could be non-zero. Thus, we set $-2(q - j) = 2i$ in \cref{eq:NonKron}. For
the current group $q = n + 1$, and therefore, $i = j - (n + 1)$. However, the upper limit of $i$ in
\cref{eq:NonKron} is $j - (M + 1)$ and also $n \le M$. All these conditions can be satisfied
together only when $n = M$. Introducing a new summation index $m = j - (M + 1)$, \cref{eq:NonKron}
simplifies to
\begin{equation}
    \notag
    \sum_{m=1}^{M}
        \alpha^{(2(m + M + 1))}_{2(M + 1), 2m} C^{(-2m)}_{2m}
            \langle \bm{p}_{2(M + 1)}, \tilde{\bm{p}}_{2M} \rangle
            \delta_{n,M}
        \text{.}
\end{equation}

Finally, the rhs of \cref{eq:genM:2q:2n:Odd} has $q = -M + s + n = n + 1$. Thus, the
inner product would not evaluate to a Kronecker delta, as that requires $q = n$, but can still be
non-zero if $q = n + 1 > M$, which is possible only for $n = M$. So, the rhs
contribution is
\begin{equation}
    \label{eq:SeqM+1:Odd:rhs}
    \beta_{2(M + 1)} \langle \bm{p}_{2(M + 1)}, \tilde{\bm{p}}_{2M} \rangle \delta_{n,M}
    \text{.}
\end{equation}

Putting all these expressions together, we obtain for the $s = M + 1$ group
\begin{equation}
\begin{split}
    \label{eq:SeqM+1}
    \sum_{m=1}^{M}
        \alpha^{(2(m + n))}_{2n, 2m} C^{(-2(m - 1))}_{2m} &= \\
            \langle \bm{p}_{2(M + 1)}, \tilde{\bm{p}}_{2M} \rangle \biggl(
                &\Bigl( \beta_{2(M + 1)}
                       - \sum_{m=1}^{M} \alpha^{(2(m + M + 1))}_{2(M + 1), 2m} C^{(-2m)}_{2m}
                    \Bigr) \delta_{n,M} \\
                &- \alpha^{(2(M + 1 + n))}_{2n, 2(M + 1)} C^{(-2M)}_{2M}
            \biggr)
\end{split}
\end{equation}
for $n = 1, \ldots, M$. The matrix of $\alpha^{(2(m + n))}_{2n,2m}$ in the lhs of
\cref{eq:SeqM+1} is the same as in \cref{eq:SeqM}.

For $M = 2$, \cref{eq:SeqM+1} reproduces \cref{eq:hstar4:p2:m2,eq:hstar6:p4:m2}, while for $M = 3$,
it reproduces \cref{eq:hstar4:p2:m3,eq:hstar6:p4:m3,eq:hstar8:p6:m3}.

\subparagraph{Groups $s > M + 1$}

By now explicit expansion of various terms in \cref{eq:genM:2q:2n:Odd} becomes too cumbersome and the
resulting expressions more complex than the original equation (compare \cref{eq:genM:2q:2n:Odd} with
\cref{eq:SeqM+1}). Thus, in the lhs of \cref{eq:genM:2q:2n:Odd} we separate the terms with
inner products evaluating to Kronecker delta functions and keep the rest as a triple sum with extra
conditions:
\begin{equation}
\begin{split}
    \label{eq:SgtM+1}
    \sum_{m=1}^{M}
        &\alpha^{(2(m + n))}_{2n, 2m} C^{(2(-M + s - m))}_{2m} =
          \beta_{2q} \langle \bm{p}_{2q}, \tilde{\bm{p}}_{2n} \rangle \\
          &- \sum_{j=2}^{s + n} \sum_{i=1}^{j-1} \sum_{l=1}^M
              \biggl.
                  \alpha^{(2j)}_{2(j - i), 2i} C^{(2(q - j))}_{2l}
                  \langle \bm{p}_{2i}, \tilde{\bm{p}}_{2l} \rangle
                  \langle \bm{p}_{2(j - i)}, \tilde{\bm{p}}_{2n} \rangle
              \biggr|_{\raisemath{1pt}{\text{per \cref{eq:SgtM+1:Cond}}}}
\end{split}
\end{equation}
for $n = 1, \ldots, M$ and $q = -M + s + n > n + 1 > 1$. The logical conditions in the bar subscript
select the second and third types of summands described in the previous subsection:
\begin{equation}
    \label{eq:SgtM+1:Cond}
    \begin{aligned}
        j - i & \le M \; \& \; i > M \\
              & \text{or} \\
        j - i & > M \text{.}
    \end{aligned}
\end{equation}
All the variables $C^{(\ldots)}_{2i}$ in those terms are already known from the
earlier groups,
although they may or may not be zeros. The matrix in the lhs of \cref{eq:SgtM+1} is the
same as in \cref{eq:SeqM,eq:SeqM+1}, so, it can be inverted only once allowing a more efficient
numerical implementation.

The presented \cref{eq:SlessM,eq:SeqM,eq:SeqM+1,eq:SgtM+1} technically allow computing the power
series of $\bm{w}$ up to any desired order, although the convergence of such series has not been
established.

\paragraph{Expressions for $L^2_{2k}$\label{sec:SmallH:L2}}

Given $\bm{w}$ and $\bm{b}$, the GS mismatch is found as $L^2 = {\lVert F \rVert}_2 - \bm{w}\tran
\bm{b}$ (see \cref{G1:eq:L2:Quad}). Therefore, the power series for $L^2$ is a product of power series for
$\bm{w}$ (\cref{eq:decw:h0:Odd}, sum over the index $m$) and $\bm{b}$ (\cref{eq:decb:h0:Odd}, sum
over the index $n$). The $h_*^0$ terms appear in the product when $m = -n$. By the limits on $n$ and
$m$, there are $M$ such instances at $-M \le m \le -1$. Taking all such terms together and using
\cref{eq:decw2m:h0:Odd}, we obtain
\begin{equation}
    \notag
    \begin{aligned}
    (\bm{w}\tran \bm{b})_0 &= \sum_{n=1}^M \beta_{2n} \bm{w}_{-2n}\tran \bm{p}_{2n} \\
                           &= \sum_{n=1}^M \beta_{2n} \sum_{j=1}^M C_{2j}^{(-2n)} \langle \bm{p}_{2n}, \tilde{\bm{p}}_{2j} \rangle
                           \text{.}
    \end{aligned}
\end{equation}
As $n, j \in [1, M]$, only the terms with $j = n$ in the second sum are non-zero, and we find
\begin{equation}
    \label{eq:l2hstar0:mM}
    L^2_0 = {\lVert F \rVert}_2 - \sum_{n=1}^M \beta_{2n} C_{2n}^{(-2n)} \text{,}
\end{equation}
with ${\lVert F \rVert}_2$ from \cref{eq:F2}.
As for the special cases of $M = 2$ and $3$, all the coefficients $C_{2n}^{(-2n)}$ used in this
expression come from the solution of \cref{eq:SeqM}.

Continuing the same approach, the $h_*^2$ terms appear in the product of two sums when $m = -n + 1$,
and there are $M + 1$ such instances at $-M \le m \le 0$. Taking all such terms together and using
\cref{eq:decw2m:h0:Odd}, we obtain
\begin{equation}
    \notag
    \begin{aligned}
    L^2_2 = - (\bm{w}\tran \bm{b})_2
        &= -\sum_{n=1}^{M + 1} \beta_{2n} \bm{w}_{-2(n - 1)}\tran \bm{p}_{2n} \\
        &= -\sum_{n=1}^{M + 1}\beta_{2n} \sum_{j=1}^M C_{2j}^{(-2(n - 1))} \langle \bm{p}_{2n}, \tilde{\bm{p}}_{2j} \rangle \\
        &= -\sum_{n=1}^{M} \beta_{2n} C_{2n}^{(-2(n - 1))} - \beta_{2(M + 1)} \sum_{j=1}^M C_{2j}^{(-2M)} \langle \bm{p}_{2(M + 1)}, \tilde{\bm{p}}_{2j} \rangle \text{.}
    \end{aligned}
\end{equation}
The second term can be simplified because $C_{2j}^{(-2M)} = 0$ for $j < M$. This yields
\begin{equation}
    \label{eq:l2hstar2:mM}
    L^2_2 = -\sum_{n=1}^{M} \beta_{2n} C_{2n}^{(-2(n - 1))} - \beta_{2(M + 1)} C_{2M}^{(-2M)} \langle \bm{p}_{2(M + 1)}, \tilde{\bm{p}}_{2M} \rangle \text{.}
\end{equation}
The coefficients $C_{2n}^{(-2(n -1))}$ are found from solving \cref{eq:SeqM+1}.

To obtain the arbitrary $h_*^{2k}$, $k > 1$, term of $\bm{w}\tran \bm{b}$, we need $m = -n + k$, and
there are $M + k$ such instances for $-M \le m \le -1 + k$. Their sum is
\begin{equation}
    \notag
    \begin{aligned}
    (\bm{w}\tran \bm{b})_{2k} &= \sum_{n=1}^{M + k} \beta_{2n} \bm{w}_{-2(n - k)}\tran \bm{p}_{2n} \\
                              &= \sum_{n=1}^{M + k} \beta_{2n} \sum_{j=1}^M C_{2j}^{(-2(n - k))} \langle \bm{p}_{2n}, \tilde{\bm{p}}_{2j} \rangle \\
                              &= \sum_{n=1}^{M} \beta_{2n} C_{2n}^{(-2(n - k))} + \sum_{n=M+1}^{M+k} \beta_{2n} \sum_{j=1}^M C_{2j}^{(-2(n - k))} \langle \bm{p}_{2n}, \tilde{\bm{p}}_{2j} \rangle \text{.}
    \end{aligned}
\end{equation}
We can further simplify the second sum by removing the terms with $C_{2j}^{(-2(n - k))}$ that are
known to be zero.
When $n = M + k$, we have $C_{2j}^{(-2M)}$,
which is zero except when $j = M$; when $n = M + k - 1$, we have
$C_{2j}^{(-2(M - 1))}$, which is zero except when $j \in [M - 1, M]$, and so on, such that we can
start the second sum over $j$ from $n - k$, resulting in
\begin{equation}
    \label{eq:l2hstar2k:mM}
    L^2_{2k} = -\sum_{n=1}^{M} \beta_{2n} C_{2n}^{(-2(n - k))} - \sum_{n=M+1}^{M+k} \beta_{2n} \sum_{j=n-k}^M C_{2j}^{(-2(n - k))} \langle \bm{p}_{2n}, \tilde{\bm{p}}_{2j} \rangle \text{.}
\end{equation}

We can also prove that all terms with negative powers of $h_*$, \textit{i.e.} $(\bm{w}\tran
\bm{b})_{2k}$ for $k < 0$, are equal to zero.\footnote{Requiring these terms to be zero is an
alternative way to derive the equations for $C_{2j}^{(-2k)}$.} This is expected because $L^2$ is a finite
quantity and cannot grow to infinity when $h_* \to 0$. See \cite{Xiourouppa:2027:PhD} for details.

\paragraph{The lowest power of \texorpdfstring{$h$}{h} in the decomposition of
weights}\label{sec:SmallH:MM}

In presenting \cref{eq:decw:h0:Odd}, we stated without proof that the lowest power of $h_*$
in the power series of $\bm{w}$ is $-2M$. Now we can apply the obtained solution to prove this
statement.

The $\bm{w}_{-2M}$ vector, which multiplies $h_*^{-2M}$ in \cref{eq:decw:h0:Odd}, is a sum of
$\tilde{\bm{p}}_{2j}$ with the coefficients $C^{(-2M)}_{2j}$ for $j = 1, \ldots, M$. According to
\cref{eq:SlessM}, these coefficients are all zeros for $j < M$: consider the terms with $m = s$ for
$s = 1, \ldots, M - 1$. The only non-zero coefficient in $\bm{w}_{-2M}$ is $C^{(-2M)}_{2M}$; it
comes from the $m = M$ term of \cref{eq:SeqM}.

Now consider what would happen if we start the power series of $\bm{w}$ with $\bm{w}_{-2(M + 1)}
h_*^{-2(M + 1)}$. Applying the same decompositions of $\bm{A}$ and $\bm{b}$, we would obtain the
same equations as \labelcref{eq:SlessM,eq:SeqM}, just with a shift in the upper index of
$C^{(\ldots)}_{2j} \Rightarrow C^{(\ldots - 2)}_{2j}$ and a shift in the rhs index,
$\beta_{2q} \Rightarrow \beta_{2q - 2}$ (see \cref{eq:genM:Awbh:Odd}). Specifically for $q = 1$, the
change would be $\beta_{2} \Rightarrow 0$, resulting in zero rhs of the index-shifted
\cref{eq:SeqM}. Thus, we would have $C^{(-2(M + 1))}_{2j} = 0$ for all $j$ and $\bm{w}_{-2(M + 1)} =
\bm{0}$.
In other words, we can start the power series of $\bm{w}$ at any $m = \mathcal{M}$ such that
$\mathcal{M} < -M$, but the first non-zero term would be $\bm{w}_{-2M} = C^{(-2M)}_{2M}
\tilde{\bm{p}}_{2M}$.

\subsection{Even case}\label{sec:SmallH:Even}

Now we turn to the complementary case when the total number of terms in the GS approximation is
even. The specific definitions of $M$, $\bm{b}$, and $\bm{A}$ for this case are given in
\cref{G1:sec:1D:Even} of \cite{Mikhin:2026:Splitting}.
Due to the proliferation of terms in \cref{G1:eq:cd:Even} as compared to \cref{G1:eq:cd:Odd}, the asymptotic analysis of the even-length
case for $h \to 0$ is more cumbersome. We again search for the solution in the form of power
series over $h_*$, as in \cref{eq:decb:h0:Odd,eq:deca:h0:Odd,eq:decw:h0:Odd}. The partial
derivatives of $b_m$ and $a_{m, k}$ with respect to $h_*$ at $h_* = 0$ are (compare with the
odd-case expressions in \cref{eq:dbdh:h0:Odd,eq:dadh:h0:Odd}):
\begin{alignat}{2}
    \notag
        &\frac{1}{(2n)!}
        \frac{\partial^{2n} b_m}{\partial h_*^{2n}} \Bigr|_{h_*=0} = \;
            &&\frac{(-1)^n \sigma^{2n}}{n! \, 2^{n - 1/2} (1 + \sigma^2)^{n + 1/2} \sqrt{\pi}}
                \Bigl[ (m + \sfrac{1}{2})^{2n} - (\sfrac{1}{2})^{2n} \Bigr] \\
    \label{eq:dbdh:h0:Even}
            &
            &&+\frac{(-1)^n}{n! \, 2^{2n + 1} \sigma \sqrt{\pi}}
                \Bigl[ 1 - m^{2n} - (m + 1)^{2n} \Bigr] \text{,} \\
    \notag
        &\frac{1}{(2n)!}
        \frac{\partial^{2n} a_{m,k}}{\partial h_*^{2n}} \Bigr|_{h_*=0} = \;
            &&\frac{(-1)^n}{n! \, 2^{2n} \sigma \sqrt{\pi}}
                \Bigl[ (k - m)^{2n} + (k + m + 1)^{2n} \Bigr. \\
    \label{eq:dadh:h0:Even}
            &
            &&\quad\quad \Bigl. - (m + 1)^{2n} - (k + 1)^{2n} - m^{2n} - k^{2n} + 1 \Bigr] \text{.}
\end{alignat}
We use a new definition for the $\bm{p}_{2n}$ vectors (compare with \cref{eq:p2n:Odd}) with
the $m^\text{th}$ elements equal to
\begin{equation}
    \label{eq:p2n:Even}
    {\left( \bm{p}_{2n} \right)}_m = (m + \sfrac{1}{2})^{2n} - (\sfrac{1}{2})^{2n} \text{.}
\end{equation}
As in the odd case, these vectors are independent for $n \le M$, but not mutually orthogonal.
Then the new complementary vectors $\tilde{\bm{p}}_{2n}$ for $n = 1, ..., M$ are defined by the
same rules as earlier: orthogonal to all $\bm{p}_{2j}$ for $j < M$ and $j \ne n$, and normalised
such that $\langle \tilde{\bm{p}}_{2n}, \bm{p}_{2n} \rangle = 1$.

Introducing the notation $\nu_m \defeq m + \sfrac{1}{2}$, $\forall m$, and assuming $n > 1$, the last term in the
square brackets in \cref{eq:dbdh:h0:Even} is represented as
\begin{align}
    \notag
        1 - &m^{2n} - (m + 1)^{2n} \\
    \notag
            &= 1 - \left( \nu_m - \sfrac{1}{2} \right)^{2n} - \left( \nu_m + \sfrac{1}{2} \right)^{2n} \\
    \notag
            &= 1 - 2 \sum_{j = 0}^n \binom{2n}{2j} \nu_m^{2j} \frac{1}{2^{2n - 2j}} \\
    \notag
            &= 1 - 2 \biggl[
                \nu_m^{2n} - \frac{1}{2^{2n}} + \frac{2}{2^{2n}}
                    + \sum_{j = 1}^{n - 1} \binom{2n}{2j}
                        \left( \nu_m^{2j} - \frac{1}{2^{2j}} + \frac{1}{2^{2j}} \right) \frac{1}{2^{2n - 2j}}
                \biggr] \\
    \notag
            &= 1 - 2 \biggl[
                \left( \bm{p}_{2n} \right)_m + \frac{2}{2^{2n}}
                    + \sum_{j = 1}^{n - 1} \binom{2n}{2j}
                        \left( \left( \bm{p}_{2j} \right)_m + \frac{1}{2^{2j}} \right) \frac{1}{2^{2n - 2j}}
                \biggr] \\
    \label{eq:dbdh:h0:Even:nu}
            &= -2 \left( \bm{p}_{2n} \right)_m
                    - 2 \sum_{j = 1}^{n - 1} \binom{2n}{2j} \frac{ \left( \bm{p}_{2j} \right)_m }{2^{2n - 2j}}
        \text{.}
\end{align}
For the last-line simplification we used the identities
\begin{equation}
    \begin{aligned}
    \notag
    \sum_{j = 0}^{n} \binom{2n}{2j}
        &= \sum_{j = 0}^{n} \binom{2n}{2j} 1^{2n - 2j} 1^{2j}
        = \frac{ (1 - 1)^{2n} + (1 + 1)^{2n} }{2} = 2^{2n - 1} \text{, and} \\
    \notag
    \sum_{j = 1}^{n - 1} \binom{2n}{2j}
        &= \sum_{j = 0}^{n} \binom{2n}{2j} - 2 = 2^{2n - 1} - 2 \text{.}
    \end{aligned}
\end{equation}
Adopting the usual interpretation that the sum is zero if the upper limit is below the lower,
expression \labelcref{eq:dbdh:h0:Even:nu} is also valid for the special case of $n = 1$, as
\begin{equation}
    \notag
    1 - m^{2n} - (m + 1)^{2n} = -2 m^2 - 2 m = -2 \left( \bm{p}_{2} \right)_m \text{.}
\end{equation}
Thus, substituting \cref{eq:dbdh:h0:Even:nu} into \cref{eq:dbdh:h0:Even}, we obtain the decomposition of
the rhs vector $\bm{b}$ as
\begin{equation}
    \label{eq:decb:h0:Even}
    \bm{b} = \sum_{n=1}^{+\infty} h_*^{2n} \sum_{j = 1}^{n} \beta^{(2j)}_{2n} \bm{p}_{2j} \text{,}
\end{equation}
where
\begin{equation}
    \label{eq:beta:h0:Even:2n}
         \beta^{(2n)}_{2n} = \frac{(-1)^n}{n! \, 2^{2n} \sigma \sqrt{\pi}} \left( \frac{2^{n + 1/2} \sigma^{2n + 1}}{(1 + \sigma^2)^{n + 1/2}} - 1 \right)
         \text{,}
\end{equation}
the same as the odd-case value from \cref{eq:beta:h0:Odd}, and
\begin{equation}
    \label{eq:beta:h0:Even:2j}
         \beta^{(2j)}_{2n} = \frac{(-1)^{n + 1}}{n! \, 2^{4n - 2j} \sigma \sqrt{\pi}} \binom{2n}{2j} \text{,}
         \quad \forall j < n
         \text{.}
\end{equation}

To simplify \cref{eq:dadh:h0:Even}, we first establish the relationship
\begin{equation}
    \notag
    \begin{aligned}
        \left( \bm{p}_{2n - 2j} \otimes \bm{p}_{2j} \right)_{(m,k)} &=
            \Bigl( \nu_m^{2n - 2j} - \sfrac{1}{2^{2n - 2j}} \Bigr) \Bigl( \nu_k^{2j} - \sfrac{1}{2^{2j}} \Bigr) \\
            &= \nu_m^{2n - 2j} \nu_k^{2j} - \nu_m^{2n - 2j} \frac{1}{2^{2j}} - \nu_k^{2j} \frac{1}{2^{2n - 2j}} + \frac{1}{2^{2n}}
        \text{.}
    \end{aligned}
\end{equation}
It is valid for $0 < j < n$ because $\bm{p}_{2k}$ is defined only for $k > 0$.

Assuming $n > 1$, we transform the square-bracket term in \cref{eq:dadh:h0:Even} as
\begin{align}
    \notag
        (k - &m)^{2n} + (k + m + 1)^{2n} - (m + 1)^{2n} - m^{2n} - (k + 1)^{2n} - k^{2n} + 1 \\
    \begin{split}
    \notag
    \begin{alignedat}[b]{2}
        &= 1 &&+ (\nu_k - \nu_m)^{2n} + (\nu_k + \nu_m)^{2n} \\
        &    &&- (\nu_m + \sfrac{1}{2})^{2n} - (\nu_m - \sfrac{1}{2})^{2n} \\
        &    &&- (\nu_k + \sfrac{1}{2})^{2n} - (\nu_k - \sfrac{1}{2})^{2n} \\
    \end{alignedat}
    \end{split} \\
    \notag
        &= 1 + 2 \sum_{j = 0}^{n} \binom{2n}{2j} \Bigl[
            \nu_m^{2n - 2j} \nu_k^{2j} -
            \nu_m^{2n - 2j} \frac{1}{2^{2j}} -
            \nu_k^{2j} \frac{1}{2^{2n - 2j}}
           \Bigr] \\
    \notag
        &= 1 + 2 \left( -\frac{2}{2^{2n}} \right)
             + 2 \sum_{j = 1}^{n - 1} \binom{2n}{2j} \Bigl[
            \left( \bm{p}_{2n - 2j} \otimes \bm{p}_{2j} \right)_{(m,k)} - \frac{1}{2^{2n}}
           \Bigr] \\
    \label{eq:dadh:h0:Even:nu}
        &= 2 \sum_{j = 1}^{n - 1} \binom{2n}{2j} \left( \bm{p}_{2n - 2j} \otimes \bm{p}_{2j} \right)_{(m,k)}
        \text{.}
\end{align}
This expression is also formally valid for the special case of $n = 1$ because
\begin{equation}
    \notag
        (k - m)^{2} + (k + m + 1)^{2} - (m + 1)^{2} - m^{2} - (k + 1)^{2} - k^{2} + 1 = 0
\end{equation}
by direct decomposition.

Substituting \cref{eq:dadh:h0:Even:nu} into \cref{eq:dadh:h0:Even}, we obtain exactly the same form
of the power series for the lhs matrix $\bm{A}$ as given by
\crefrange{eq:deca:h0:Odd}{eq:alpha:h0:Odd} for the odd case, just with a different definition of
the $\bm{p}_{2n}$ vectors, \cref{eq:p2n:Even} instead of \cref{eq:p2n:Odd}. As in the odd case, the
power series of $\bm{b}$ starts at $h_*^2$, while the series of $\bm{A}$ only at $h_*^4$. Therefore,
we expect the mixand weights to grow infinitely by modulo as $h_* \to 0$, and some of them to become
negative to comply with the normalisation condition.

Together, the presented derivation yields
\begin{equation}
    \label{eq:genM:Awbh:Even}
    \sum_{j=2}^{q + M} \bm{A}_{2j} \bm{w}_{2(q - j)} =
        \sum_{j = 1}^{q} \beta^{(2j)}_{2q} \bm{p}_{2j} [ q \geq 1 ] \text{,}
\end{equation}
and, repeating the transformation of the lhs from \cref{eq:genM:2q:Odd}, the even-case variant of
the $2q / 2n$ equations (compare with \cref{eq:genM:2q:2n:Odd}):
\begin{multline}
    \label{eq:genM:2q:2n:Even}
        \sum_{j=2}^{q + M} \sum_{i=1}^{j-1} \sum_{l=1}^M
            \alpha^{(2j)}_{2(j - i), 2i} C^{(2(q - j))}_{2l}
            \langle \bm{p}_{2i}, \tilde{\bm{p}}_{2l} \rangle
            \langle \bm{p}_{2(j - i)}, \tilde{\bm{p}}_{2n} \rangle \\
        = \sum_{j = 1}^{q} \beta^{(2j)}_{2q} \langle \bm{p}_{2j}, \tilde{\bm{p}}_{2n} \rangle [ q \geq 1 ]
        \text{.}
\end{multline}
Solution then proceeds along the same route as in \cref{sec:SmallH:GenM} for the odd case. We
sequentially take the groups of equations \labelcref{eq:genM:2q:2n:Even} with $q = -M + s + n$ for
$s = 1, 2, \ldots$ and $n$ running from $1$ to $\min{(s,M)}$.

For $s < M$, we have $q = n + (s - M) < n < M$, so that the rhs of \cref{eq:genM:2q:2n:Even} is
zero. Therefore, we end up with linear systems of equations for $C^{(2j)}_{2m}$ that are formally
identical\footnote{Equations look exactly the same, just with different definitions for
$\bm{A}_{2k}$, $\bm{p}_{2k}$, and $\tilde{\bm{p}}_{2k}$.} to the odd case \cref{eq:Seq1,eq:SlessM}.
The similarity continues for $s = M$ and $q = n$: the rhs of \cref{eq:genM:2q:2n:Even} is now
$\beta^{(2n)}_{2n}$, which by \cref{eq:beta:h0:Even:2n} is the same as $\beta_{2n}$ in
\cref{eq:SeqM}.

The first substantial difference occurs in the $s = M + 1$ group: for $q = n + 1$ the rhs of
\cref{eq:genM:2q:2n:Even} may have two non-zero terms at $j = n$ and $j = n + 1$, the latter
possible only when $n = M$. Thus, instead of \cref{eq:SeqM+1:Odd:rhs}, we get
\begin{equation}
    \label{eq:SeqM+1:Even:rhs}
    \beta^{(2n)}_{2(n + 1)} + \beta^{(2(M + 1))}_{2(M + 1)} \langle \bm{p}_{2(M + 1)}, \tilde{\bm{p}}_{2M} \rangle \delta_{n,M}
\end{equation}
and end up with a system of linear equations like \cref{eq:SeqM+1}, but with an extra term
$\beta^{(2n)}_{2(n + 1)}$ in the rhs.

A similar transformation occurs for the groups $s > M + 1$, which result in equations similar to
\cref{eq:SgtM+1}, but with the first term on the rhs replaced by
\begin{equation}
    \label{eq:SgtM+1:Even:rhs}
    \sum_{j = 1}^{q} \beta^{(2j)}_{2q} \langle \bm{p}_{2j}, \tilde{\bm{p}}_{2n} \rangle \text{.}
\end{equation}
With these relatively minor adjustments, we can find the power series of $\bm{w}$
up to any desired order in $h_*$, as we did in the odd case. Then the power series of $L^2$ are
obtained as in Section~\ref{sec:SmallH:L2}.
However, the presence of all $\bm{p}_{2j}$, $j \le n$, in the
decomposition of $\bm{b}_{2n}$ by \cref{eq:decb:h0:Even} precludes the simplification of the series
products achieved in \cref{eq:l2hstar0:mM,eq:l2hstar2:mM}. For this reason we stop at the generic expression
\begin{equation}
    \label{eq:l2hstar:2k}
    L^2_{2k} = {\lVert F \rVert}_2 \delta_{0,k}
        - \sum_{n=1}^{M + k} \bm{w}_{-2(n - k)}\tran
            \biggl( \sum_{j=1}^n \beta^{(2j)}_{2n} \bm{p}_{2j} \biggr)
    \text{,} \quad \forall k \ge 0 \text{.}
\end{equation}
The proof of the lowest power of $h_*$ in \cref{eq:decw:h0:Odd} does not change as compared to the
odd case.

\section{Derivation of \texorpdfstring{$L^2$}{L²} for large
\texorpdfstring{$M$}{M}}\label{sec:AsympM:L2}

\renewcommand{\theequation}{B.\arabic{equation}}
\setcounter{equation}{0}

\subsection{Odd case}\label{sec:AsympM:L2:Odd}

We use the cosine double-angle formula to transform \cref{G1:eq:L2:Odd:AsympM:1} into:
\begin{equation}
    \label{eq:L2:Odd:AsympM:1}
        \lim_{M \to \infty} L^2 = \frac{1}{\pi} e^{-4 \pi^2 \sigma_c^2}
            \int_{-\infty}^{+\infty} \left[e^{-x^2} + e^{-x^2} \cos{\! \left(4 \pi \frac{1 - \sigma_c^2}{h^2} x \right)} \right] \! dx \text{.}
\end{equation}
The first term in square brackets is the Gaussian integral, while the second term is known from
\cite[Eq.~7.4.6]{Abramowitz:1965:Handbook}; to account for integration from $-\infty$, we double the
result, as the integrand is an even function of $x$. Taken together, we obtain
\begin{equation}
    \label{eq:L2:Odd:AsympM}
        \lim_{M \to \infty} L^2 \approx \frac{1}{\sqrt{\pi}} e^{-4 \pi^2 \sigma_c^2} \left( 1 + e^{-\frac{4 \pi^2 (1 - \sigma^2)^2}{h^2}} \right)
        \! \text{.}
\end{equation}

\subsection{Even case}\label{sec:AsympM:L2:Even}

We start with \cref{G1:eq:L2:Even:AsympM:1} of \cite{Mikhin:2026:Splitting}, which is nearly
identical to its odd-case equivalent in \cref{G1:eq:L2:Odd:AsympM:1}. The only difference is that
the $m_c$ factor in \cref{G1:eq:mc:Even} has an extra constant of $\sfrac{-1}{2}$ as compared to
\cref{G1:eq:statec:Odd:2}. However, once we apply the cosine double-angle formula and substitute
$m_c$, the cosine argument becomes
\begin{equation}
\notag
4 \pi \frac{1 - \sigma_c^2}{h^2} x - 2 \pi \text{,}
\end{equation}
\textit{i.e.} simply shifted by a period as compared to \cref{eq:L2:Odd:AsympM:1}. The shift does not change
the function value, and therefore, integration using \cite[Eq.~7.4.6]{Abramowitz:1965:Handbook}
again produces \cref{eq:L2:Odd:AsympM}.

\bibliographystyle{alpha}
\bibliography{gs-1D-bonus}

\end{document}